\documentclass[11pt]{article}
\usepackage{amsmath}
\usepackage{cite} 
\usepackage{graphics} 
  \usepackage{epsfig} 

\usepackage[pagewise]{lineno} 

\usepackage{amsmath}
\usepackage{paralist}
\usepackage{cite}  
\usepackage{graphics} 
  \usepackage{epsfig} 
\usepackage[colorlinks=true]{hyperref}
\hypersetup{urlcolor=blue, citecolor=red}

\usepackage{amssymb}
\usepackage{upgreek} 

\usepackage{tikz}  
\usetikzlibrary{decorations.pathmorphing}
\usetikzlibrary{patterns}
\usepackage{xcolor} 
\usepackage{pgfplots}
\usetikzlibrary{calc}

	\newcommand{\COMMENT}[1]{}
	\newcommand{\DELETE}[1]{}

        \newcommand{\REM}[1]{\marginpar{\bfseries\tiny{\color{blue}}}}

\usepackage{mathrsfs}   

\newtheorem{proposition}{Proposition}

\newtheorem{remark}{Remark}

\newcommand{\Item}[2]{\parbox[t]{.05\textwidth}{#1}\hfill%
      \parbox[t]{.95\textwidth}{#2}\vspace*{.8mm}}
\newcommand{\eq}[1]{(\ref{#1})}
\newcommand{\bbD}{\mathbb D}
\newcommand{\bbH}{\mathbb H}
\newcommand{\bbI}{\mathbb I}
\newcommand{\In}{\!\in\!}
\newcommand{\pl}{\partial}
\renewcommand{\d}{\mathrm d}  
\newcommand\DT[1]{\mathchoice
                 {{\buildrel{\hspace*{.1em}\text{\LARGE.}}\over{#1}}}
                 {{\buildrel{\hspace*{.1em}\text{\Large.}}\over{#1}}}
                 {{\buildrel{\hspace*{.1em}\text{\large.}}\over{#1}}}
                 {{\buildrel{\hspace*{.1em}\text{\large.}}\over{#1}}}}
\newcommand\DDT[1]{\mathchoice
   {{\buildrel{\hspace*{.13em}\text{\LARGE.\hspace*{-.13em}.}}\over{#1}}}
   {{\buildrel{\hspace*{.1em}\text{\Large.\hspace*{-.1em}.}}\over{#1}}}
   {{\buildrel{\hspace*{.1em}\text{\large.\hspace*{-.1em}.}}\over{#1}}}
   {{\buildrel{\hspace*{.1em}\text{\large.\hspace*{-.1em}.}}\over{#1}}}}
\newcommand\DDDT[1]{\mathchoice
   {{\buildrel{\hspace*{.1em}\text{\LARGE.\hspace*{-.1em}\LARGE.\hspace*{-.1em}.}}\over{#1}}}
   {{\buildrel{\hspace*{.1em}\text{\Large.\hspace*{-.1em}\Large.\hspace*{-.1em}.}}\over{#1}}}
   {{\buildrel{\hspace*{.1em}\text{\large.\hspace*{-.1em}\large.\hspace*{-.1em}.}}\over{#1}}}
   {{\buildrel{\hspace*{.1em}\text{\large.\hspace*{-.1em}\large.\hspace*{-.1em}.}}\over{#1}}}}
\def\R{\mathbb R}
\def\N{\mathbb N}
\newcommand{\lineunder}[2]{\LU{\begin{array}[t]{c}\underbrace{#1}\vspace*{.5em}\end{array}}{\mbox{\footnotesize\rm #2}}}
\newcommand{\LU}[2]{\begin{array}[t]{c}#1\vspace*{-1em}\\_{#2}\end{array}}
\newcommand{\linesunder}[3]{\LSU{\begin{array}[t]{c}\underbrace{#1}\vspace*{.5em}\end{array}}{\mbox{\footnotesize\rm #2}}{\mbox{\footnotesize\rm #3}}}
\newcommand{\LSU}[3]{\begin{array}[t]{c}#1\vspace*{-1em}\\_{#2}\vspace*{-.5em}\\_{#3}\end{array}}
\newcommand{\barOmega}{\,\overline{\!\Omega\!}\,}
\newcommand{\divS}{\mathrm{div}_{\scriptscriptstyle\textrm{\hspace*{-.1em}S}}^{}}
\newcommand{\nablaS}{\nabla_{\scriptscriptstyle\textrm{\hspace*{-.3em}S}}^{}}
\newcommand{\Vdots}{\mbox{\large{:}}\hspace*{-.3em}^{^{\mbox{\large{.}}}}}

\allowdisplaybreaks

\begin{document}
\begin{sloppypar}


  \hspace*{2em}{\Large\bf From quasi-incompressible to semi-compressible fluids}\footnote{This research has been partially supported from the grants
19-04956S
(especially regarding the focus on the dynamic and nonlinear behaviour)
of Czech Science Foundation,
and the FWF/CSF project 19-29646L
(especially regarding the focuse on the large strains in materials science),
 and also from the 
institutional support RVO:\,61388998 (\v CR)
}


  \bigskip
  
  \centerline{\scshape Tom\'a\v{s} Roub\'\i\v cek}

  \medskip
\begin{center}
{\footnotesize
Mathematical Institute, Charles University,\\[-.0em]
Sokolovsk\'a 83, CZ-186~75~Praha~8, Czech Republic 
\\[-.2em]and\\[-.2em]
Institute of Thermomechanics, \,Czech Acad.\ Sci.,\\[-.2em]
Dolej\v skova 5, CZ-18200 Praha 8, Czech Republic 
}
\end{center}

\bigskip

\begin{quote}{\normalfont\fontsize{8}{10}\selectfont
    {\bfseries Abstract.}
   A new concept of semi-compressible fluids is introduced for slightly
    compressible visco-elastic fluids (typically rather liquids than gasses)
    where mass density variations are negligible in some sense, while
    being directly controlled by pressure which is very small in comparison
    with the elastic bulk modulus. The physically consistent fully Eulerian
    models with specific dispersion of pressure-wave speed are devised.
 This contrasts to the so-called quasi-incompressible fluids which
 are described not physically consistently and, in fact,
 only approximate ideally incompressible ones in the limit.
    After surveying and modifying models for the quasi-incompressible fluids,
    we eventually devise some fully convective models complying with
    energy conservation and capturing phenomena as pressure-wave propagation
    with wave-length (and possibly also pressure) dependent velocity.
    \par


}

\end{quote}

\medskip

\def\wt{\widetilde}
\def\vv{v}

\section{Introduction}\label{sect-intro}
It is well known that propagation of pressure waves
    (briefly P-waves, also referred
    to as longitudinal waves or, if waves are in acoustic frequencies, as sound)
cannot be modeled by ideally incompressible fluid models described
by incompressible Navier-Stokes equations. On the other hands, fully
compressible fluid models which counts with the mass density variations
are analytically and computationally very complicated (with e.g.\ only
weak-strong uniqueness known), which is usually not necessary in liquids.
For example, water density increases only by 5\,\% even when so extremely
compressed as in the 11\,km deep Mariana Trench in the Pacific ocean.

This paper focuses on modelling of liquids where mass density variation
can be neglected in a certain way or, in another way, directly controlled by
pressure which is very small in comparison with elastic bulk modulus.
This elastic modulus is very large (in comparison with occuring
pressure) but not $+\infty$. The goal is to allow propagation of P-waves
while not to drive the model too far from the incompressible Navier-Stokes one. 

In fluids, typically the rheological models are very different in the
deviatoric part and in the spherical (volumetric) part of the strain tensor.
The usual Newtonian ideally incompressible
fluids use the linear Stokes rheology in the deviatoric part and while the
spherical part is ideally rigid. The essence of {\it viscoelastic fluids} is
involvement of some elasticity into purely viscous/rigid rheology.
In their simplest variant, viscoelastic fluids
use the Kelvin-Voigt rheology in the spherical part and/or
the Maxwell in the deviatoric part. Here we will be focused on the
former option.

The viscoelastic fluids with a constant viscosity which
flow like a liquid but behave like an elastic solid when stretched out
have been invented by D.V.\,Boger \cite{Bog77HECV} in the context of polymer
solutions, cf.\ \cite{Jame09BF} for a historical survey. This primarily
concerns the shear part \COMMENT{???} but we will use this idea for 
the spherical (volumetric) part \COMMENT{ISN'T IT ALSO LIKE elastomers??}
where it is even more natural in most of fluids. 
The spherical elasticity response (i.e.\ compressibility) is characterized by
the bulk elastic modulus $K$, while the shear modulus is zero in
fluids so that shear waves cannot propagate (in contrast to solids).
Together with the mass density $\varrho$, this modulus determines
the speed of the P-waves as $\sqrt{K/\varrho}$. E.g.\ water 
has approximately (cf.\ Remarks~\ref{rem-K=K(pi)} and \ref{rem-K=K(pi)++}
below) $K=2.15\,$GPa and is thus much more compressible comparing e.g.\ with
steel which has $K\sim100-200\,$GPa.

One of a prominent application is in geophysics when modelling propagation of
seismic P-waves through fluidic parts of the Earth, i.e.\ outer core
(composed with molten iron and nickel) 
and oceans (composed from water with solution of salts).
In mathematical literature, viscoelastic fluids with elasticity
in the spherical part (called quasi-incompressible or 
sometimes also quasi-compressible)
are often considered
only as a regularization (or a numerical ``pressure stabilization'')
of incompressible fluids, cf.\ e.g.\ \cite{DonMar06DAAC,DonSpi11WSNS,Osko73SPQL} or, even
with twisted energetics (due to omitting Temam's force
$\frac\varrho2({\rm div}\,\vv)\,\vv$ in \eq{NS}
below)
also e.g.\ in \cite{Chor67NMS,Chor68NSNS} or \cite[Sec.\,4.2]{Proh97PQCM},
referring (although not explicitly) to the ignored mass-density
variation or partly ignored convection, cf.\ the discussion in
Sect.~\ref{sec-convective} below. Sometimes, quasi-incompressibility may refer
to mixtures and then it means that density is a function of concentration
only, cf.\ e.g.\ \cite{FeLuMa16PDEA,LowTru98QICH}.

The goal of this article is to modify the quasi-incompressible models devised
as only approximate models for the ideally incompressible one towards
physically sensible models. We also focus on identifying dispersive character
of the specific models,
which seems rather ignored in literature.
The dispersion
of velocities of P-waves is mainly manifested in solids but it is
reported in fluids too, e.g.\ in the Earth's core \cite{Dens52LWEC,Gute58WVEC}.

In particular, 
the devised models to which we
will focus will exhibit some of (or all) the following attributes:
\begin{enumerate}[(a)]
\item 
  propagation and dispersion of P-waves
  are controlled in a certain way,
\vspace*{-.6em}\item the energy balance is preserved at least formally,
  but in some models even rigorously,
\vspace*{-.6em}\item the pressure is well defined in a reasonable sense also on
  the boundary,
\vspace*{-.6em}\item the equations are consistently written in Eulerian coordinates
  (i.e.\ the model is fully convective),
\vspace*{-.6em}\item in some models, uniqueness of weak solutions holds even in the physically
  relevant 3-dimensional cases.
  \end{enumerate}
We begin with presentation of a rather simple model satisfying (a)-(b)
in Sect.\,\ref{sect-simple}. In Sect.\,\ref{sec-dispersive}
and \ref{sec-dispersive2}, we will improve it to satisfy also (c)
by including some gradient terms with a conservative or a dissipative
character, respectively.
Eventually, in Sect.\,\ref{sec-convective}, we arrive to fully convective
models which satisfy (a)--(d) and, in the case of a multipolar variant,
also (e). We will call such new fluidic models satisfying (a)--(d)
as {\it semi-compressible fluids}.


\section{Simple quasi-incompressible models}\label{sect-simple}
The ultimate departure point is the fully
compressible Navier-Stokes system which can be rigorously derived
from Boltzmann kinetic theory. 
In term of the velocity field $\vv$ and the density $\rho$, this system is
\begin{subequations}\label{compress-NS}
  \begin{align}\label{compress-NS1}
    &\DT\rho+
\vv\cdot\nabla\rho+\rho\,{\rm div}\,\vv
    =0,
    \\&\label{compress-NS2}
    \rho\DT \vv+\rho(\vv{\cdot}\nabla)\vv-{\rm div}
    \,\sigma
    =\rho g
    \\&\label{compress-NS3}
    \text{with the stress }\ \
    \sigma=\bbD e(\vv)
-\pi(\rho)\bbI,
\end{align}\end{subequations}
where
\begin{align}\label{Be(v)}
  \bbD e(\vv)= K_\text{\sc v}^{}\,{\rm sph}\,e(\vv)+G_\text{\sc v}^{}\,{\rm dev}\,e(\vv)
\end{align}
with $K_\text{\sc v}>0$ and $G_\text{\sc v}>0$ a bulk and a shear viscosity coefficients and 
with the small-strain rate
$e=e(\vv):=\frac12\nabla\vv^\top\!+\frac12\nabla\vv$
and with ${\rm sph}\,e=({\rm tr}\,e)\bbI/d$ and
${\rm dev}\,e=e-({\rm tr}\,e)\bbI/d$ are the spherical (volumetric)
and the deviatoric parts of $e$, respectively; here ``\,tr\,'' denotes
the trace of a matrix.
We can also write
$\bbD_{ijkl}=
K_\text{\sc v}^{}\delta_{ij}\delta_{kl}+G_\text{\sc v}^{}(\delta_{ik}\delta_{jl}
+\delta_{il}\delta_{jk}-2\delta_{ij}\delta_{kl}/d)$ with ``$\,\delta\,$'' denoting the
Kronecker symbol. The pressure $\pi$ is assumed to
be a function of the density $\rho$.

The quasi-incompressible model is usually derived by considering
a small-perturbation ansatz ($\varrho>0$ being a fixed constant -- i.e.\
a homogeneous fluid):
 \begin{align}\label{ansatz}
\rho=\varrho\Big(1+\frac\pi K\Big)
 \end{align}
 where $\pi$ is the pressure and $K$ the elastic bulk modulus,
 both with the physical dimension Pa=J/m$^3$.
 Let us remark that $1/K$ is called compressibility. 
 The mass conservation \eq{compress-NS1}  written in the
 form $\DT\rho+{\rm div}(\rho\vv)=0$  then results to
$\DT\pi/K+{\rm div}((1{+}\frac\pi K\big)\vv)=0$
 and \eq{compress-NS3}
 results to $\varrho\DT \vv+\varrho(\vv{\cdot}\nabla)\vv
 -{\rm div}\,\sigma/(1{+}\frac\pi K)
    =\varrho g$
 while \eq{compress-NS3} gives
 $\sigma=\bbD e(\vv)
 -\pi\bbI$. Ignoring now the mass density variation, i.e.\
 assuming $|\pi|\ll K$ and forgetting the $\pi/K$-terms,
 we obtain $\DT\pi+K\,{\rm div}\,\vv=0$
 and, from \eq{compress-NS2},
 $\varrho\DT \vv+\varrho(\vv{\cdot}\nabla)\vv-{\rm div}\,\sigma
 =\varrho g$. For brevity, in what follows, we will write
  just  $g$ instead $\varrho g$.

 This approximation of the incompressibility constraint,
 used already e.g.\ in \cite{Chor67NMS,Chor68NSNS},
 however twists the energy balance, which is desired
 to be improved by a suitable modification of the approximate system.
 In term of the velocity field $\vv$ and the pressure $\pi$ with
 $\varrho>0$ constant, the governing
system in its simple variant is considered as
\begin{subequations}\label{NS}
  \begin{align}
    &\varrho\DT \vv+\varrho(\vv{\cdot}\nabla)\vv-{\rm div}
\big(\bbD e(\vv)-\pi\bbI\big)
=g-\frac\varrho2({\rm div}\,\vv)\,\vv
\label{NS-a}
\,,
    \\\label{NS-b}
    &\DT\pi+K\,{\rm div}\,\vv\,=0\,.
\end{align}\end{subequations}
Here $\bbI\in\R^{d\times d}$ denotes the unit matrix. Beside, $g$ is a given
bulk force, e.g.\ the gravitational force, and $K$ is the mentioned bulk
modulus. The force
$-\frac\varrho2({\rm div}\,\vv)\,\vv$ in \eq{NS-a} is related with the
convective term $\varrho(\vv{\cdot}\nabla)\vv$ in the only quasi- (but not fully)
incompressible model and was invented (for $\varrho=1$) by R.\,Temam
\cite{Tema69ASEN}
(cf.\ also \cite[Ch.\,III,\,Sect.\,8]{Tema77NSE}) and used also
e.g.\ in \cite{BerSpi17CSWS,DonMar06DAAC,DonSpi11WSNS,ChLaML19AVSN}.
A physical meaning
in the energy balance of this term is the variation of the kinetic energy
$\frac\varrho2|\vv|^2$ within the volume changes ${\rm div}\,\vv$ and is 
needed to gain the correct energetics, cf.\ \eqref{convective-tested} below.
Simultaneously, this forces vanishes in the incompressible limit.
Recently, some justification by pullback/pushforward geometrical arguments
has been provided in \cite{Toma??ITST}, based on \cite{podi97II}.

We will consider the system \eq{NS} on 
a bounded domain $\Omega\subset\R^d$ evolving in time up to $T>0$ a finite
fixed time horizon. We abbreviate $I:=[0,T]$.

The system \eqref{NS} is to be completed by the boundary conditions, say
a combination of the Dirichlet and the Navier conditions:
  \begin{align}\label{BC}
    \vv_{\rm n}=0\qquad\text{and}\qquad(\sigma\vec{n})_{\rm t}
    +\kappa\vv_{\rm t}=0,\ \ 
  \end{align}
  where $\vv_{\rm n}=\vv\cdot\vec{n}$ is the normal velocity (as a scalar)
  with $\vec{n}=\vec{n}(x)$ denoting the unit outward normal to $\Gamma$,
  while $(\cdot)_{\rm t}$ denotes the tangential component of a vector,
  i.e.\ e.g.\
  $\vv_{\rm t}=\vv-\vv_{\rm n}\vec{n}$ is the tangential velocity (a vector).
  Furthermore,
  $\sigma$ in \eq{BC} is the stress $\bbD e(\vv)-\pi\bbI$ occurring
  \eqref{NS-a}, and with $\kappa$ the ``friction'' coefficient on the boundary
  acting in the Navier boundary condition \eqref{BC}.
Eventually, considering the initial-value problem, we also prescribe the
  initial conditions
\begin{align}\label{IC}
  \vv(0)=\vv_0\qquad\text{and}\qquad \pi(0)=\pi_0
  \,.
\end{align}

The energetics of the model \eqref{NS}--\eqref{BC}
is revealed by testing \eqref{NS-a} by $\vv$.
We use \eqref{NS-b} together with the
calculus for the convective term 
\begin{align}\nonumber
  \int_\Omega\varrho(\wt\vv{\cdot}\nabla)\vv\cdot \vv\,\d x
 & =\int_\Gamma\varrho|\vv|^2(\wt\vv\cdot\vec{n})\,\d S
  - \int_\Omega\varrho\vv\nabla(\wt\vv\otimes\vv)\,\d x
  \\&\nonumber=\int_\Gamma\varrho|\vv|^2(\wt\vv\cdot\vec{n})\,\d S
  - \int_\Omega\varrho|\vv|^2{\rm div}\,\wt\vv
  +\varrho\wt\vv{\cdot}\nabla\vv\cdot\vv\,\d x
  \\&=\int_\Gamma\frac\varrho2|\vv|^2(\wt\vv\cdot\vec{n})\,\d S
  -\int_\Omega\frac\varrho2|\vv|^2({\rm div}\,\wt\vv)\,\d x
  \label{convective-tested}\end{align}
which, when used for $\wt\vv=\vv$, gives rise to the bulk force
$\frac12\varrho({\rm div}\,\vv)\vv$. Let us note that, in
\eq{convective-tested}, we exploited that $\varrho$ is constant.
Moreover, we use the calculus of the pressure gradient tested by $\vv$:
  \begin{align}\label{calculus-pi}
  \int_\Omega\nabla\pi{\cdot}\vv\,\d x-\int_\Gamma\!\pi\,(\vv{\cdot}\vec{n})\,\d S=
  -\int_\Omega\!\pi\,{\rm div}\,\vv\,\d x=
    \int_\Omega\pi\frac{\DT\pi}K\,\d x
    =\frac{\d}{\d t}\int_\Omega\frac{\pi^2}{2K}\,\d x\,.
  \end{align}
The mentioned test thus gives (at least formally) the energy balance
\begin{align}\nonumber
&\int_\Omega\!\!\!\!\linesunder{\frac\varrho2|\vv(t)|^2}{kinetic}{energy}\!\!\!\!\!+
  \!\!\!\!\!\linesunder{
    \frac1{2K}\pi(t)^2
  }{stored}{energy}\!\!\!\!\d x
  +\int_0^t\!\!\int_\Omega\!\!\!\!\!\!\linesunder{\bbD e(\vv){:}e(\vv)
  }{dissipation}{by viscosity}
  \!\!\!\d x\d t
  \\&\qquad\qquad
  +\int_0^t\!\!\int_\Gamma\!\!\!\!\!\!\!\!\!\!\!\!\!\!\!\!\linesunder{\kappa|v_{\rm t}|^2}{dissipation by a}{boundary ``friction''}\!\!\!\!\!\!\!\!\!\!\!\!\!\!\!\!\,\d S\d t
  =\int_\Omega\frac\varrho2|\vv_0|^2+\frac1{2K}\pi_0^2
  \,\d x+\int_0^t\!\!\int_\Omega\!\!\!\!\!\!\linesunder{g\cdot \vv}{power}{of loading}\!\!\!\!\!\!\d x\d t
  \,.
  \label{energy}\end{align}

We will use the standard notation concerning the Lebesgue and the Sobolev
spaces, namely $L^p(\varOmega;\R^n)$ for Lebesgue measurable functions $\varOmega\to\R^n$ whose Euclidean norm is integrable with $p$-power, and
$W^{k,p}(\varOmega;\R^n)$ for functions from $L^p(\varOmega;\R^n)$ whose
all derivative up to the order $k$ have their Euclidean norm integrable with
$p$-power. We also write briefly $H^k=W^{k,2}$. Moreover, for a Banach space $X$
and for $I=[0,T]$,
we will use the notation $L^p(I;X)$ for the Bochner space of Bochner
measurable functions $I\to X$ whose norm is in $L^p(I)$, 
and $H^1(I;X)$ for functions $I\to X$ whose distributional derivative
is in $L^2(I;X)$. Furthermore, $C_{\rm w}(I;X)$ will denote the
Banach space of weakly continuous functions $I\to X$.

Assuming $g\in L^1(I;L^2(\Omega;\R^d))+L^2(I;L^{6/5}(\Omega;\R^d))$,
$\vv_0\in  L^2(\Omega;\R^d)$, and $\pi_0\in L^2(\Omega)$,
and based on this energy balance \eqref{energy} as an inequality, by
usual approximation methods and subsequent limit passage
one can easily see existence of weak solutions
$(\vv,\pi)$ with
\begin{align}\label{est}
\vv\in
L^\infty(I;L^2(\Omega;\R^d))\:\cap\:L^2(I;H^1(\Omega;\R^d))
\quad\text{ and }\quad\pi\in L^\infty(I;L^2(\Omega))\,.
\end{align}
{}The qualification of $g$ uses
$L^{6/5}(\Omega;\R^d)\subset H^1(\Omega;\R^d)^*$ 
and the a-priori estimates \eq{est} are obtained from \eq{energy} by
using Gronwall inequality (cf.\ \cite[Lemma 8.26 with $p{=}2$ and
$q{=}\infty$]{Roub13NPDE}  for details)
combined with the Korn inequality. Moreover, 
by comparison  such weak solution  also satisfies
\begin{align}\label{est+}
  \varrho\DT \vv\in L^2(I;H^1(\Omega;\R^d)^*)+L^{5/4}(I{\times}\Omega;\R^d)
  \quad\text{ and }\quad
\DT\pi=-{\rm div}\,\vv\in L^2(I{\times}\Omega)\,.
\end{align}
The former estimate in \eq{est+}
needs $d\le3$ and uses the embedding $H^1(\Omega)\subset L^6(\Omega)$
  and the interpolation $L^\infty(I;L^2(\Omega;\R^d))\cap L^2(I;L^6(\Omega;\R^d))\subset
  L^{10/3}(I{\times}\Omega;\R^d)$ so that both $\varrho(\vv{\cdot}\nabla)\vv$ and
  $\frac12\varrho({\rm div}\,\vv)\vv$ belong to $L^{5/4}(I{\times}\Omega;\R^d)$.


Let us remark that the estimate \eq{est} still does not yield
enough integrability of the pressure to allow for uniqueness in the
physically relevant 3-dimensional case, while only 2-dimensional case
works by estimation like in
Proposition~\ref{prop-1}(iii)
below. Therefore, it is
not much improvement of the incompressible model where the uniqueness
problem has been recognized extremely difficult by involving into seven
millennium problems by the Clay institute.

\begin{remark}[{\sl Energy conservation alternatively}]\label{rem-conserv}
  \upshape
  Using the calculus
  $(\vv{\cdot}\nabla)\vv={\rm div}(\vv\otimes\vv)-({\rm div}\,\vv)\vv$,
the force balance \eq{NS-a} can also be written as
\begin{align}\label{variant-1}
\varrho\DT\vv-{\rm div}(\bbD e(\vv)-\varrho\vv\otimes\vv-\pi\bbI)
=g+\frac\varrho2({\rm div}\,\vv)\,\vv,
\end{align}
as pointed out in \cite[Rem.1.2]{BerSpi17CSWS}. Formally,
the energetics could be balanced also by other ways, e.g.\ by
the additional pressure instead of the force
\begin{align}\label{variant-2}
  &\varrho\DT \vv
  -{\rm div}\Big(\bbD e(\vv)
  -\varrho\vv\otimes\vv+\Big(\frac\varrho2|\vv|^2-\pi\Big)\bbI\Big)
  =g+\varrho({\rm div}\,\vv)\,\vv\,.
\end{align}
One could also combine
\eq{variant-1} and \eq{variant-2} with suitable weights 2 and $-1$ to
 obtain 
\begin{align}\label{variant-3}
\varrho\DT \vv-{\rm div}\Big(\bbD e(\vv)-\varrho\vv\otimes\vv-
\Big(\frac\varrho2|\vv|^2+\pi\Big)\bbI\Big)=g\,.
\end{align}
Although the desired energy is conserved in \eq{variant-2} and
\eq{variant-3} and the latter case even uses the enhanced pressure
$\frac\varrho2|\vv|^2+\pi$ like in \eq{NS-modif+++} below,
both these variants are physically not justified
 and seems not relevant because \eq{variant-2} 
is,
in fact, equivalent by replacing the force
$-\frac\varrho2({\rm div}\,\vv)\,\vv$ in \eq{NS-a} by
$\varrho(\nabla\vv)^\top\vv$  which does not vanish for
incompressible limit. 
\end{remark}

\section{Towards 
  dispersion through conservative gradient terms}\label{sec-dispersive}

The P-waves in real media usually does not propagate with a
constant speed, but the speed may depend on their frequency (or
equivalently on their wave length). This is referred to as a
dispersion and some dispersion is mostly a real effect 
thus desirable to reflect it in the model in some ways.
If the speed increases (or decreases) with the wave length, we speak
about normal (or anomalous) dispersion, respectively.
Already due to the attenuation by viscous Kelvin-Voigt-like
rheology in Sect.~\ref{sect-simple} controlled by $K_\text{\sc v}$
causes a normal dispersion,
cf.\ \cite[Rem.\,6.4.13]{KruRou19MMCM}, but usual modelling
ansatz is that $K_\text{\sc v}>0$ is rather very small to allow propagation
of waves which have also high frequencies without any essential
attenuation.  For some fluidic applications where P-waves should
propagate for long distances (say thousands of kilometers as in Earth's outer
core or oceans),  it is thus desirable to have another mechanisms
to control dispersion at disposal, not directly causing atenuation.

A certain way to realize it is a Sobolev-type regularization of the
incompressibility constraint
\begin{align}\label{Sobolev}
\frac1K\DT\pi+{\rm div}\,\vv=\frac1H\Delta\DT\pi
\end{align}
where $H$ is an elastic bulk ``hyper-modulus'' having the physical
dimension Pa/m$^2$=J/m$^5$.
The parabolic equation of the type ``$\,\DT\pi{-}\Delta\DT\pi=f\,$'' is
sometimes referred to as a Sobolev equation, which is why 
\eq{Sobolev} is called a Sobolev-type regularization here.
The term $\Delta\DT\pi$ in \eq{Sobolev} has been considered in
\cite[Sect.4.3]{Proh97PQCM} but without the $K$-term for merely numerical
purposes.

The Navier boundary conditions \eq{BC} can be now generalized by
considering the pressure-dependent friction coefficient $\kappa=\kappa(\pi)$.
Moreover, some boundary conditions are now needed for $\DT\pi$. One can consider
e.g.\ the elastical wall (boundary), i.e.\ Newton-Robin-type condition
$\nabla\DT\pi\cdot\vec{n}/H+\DT\pi/c=0$.
Equivalently, we can formulate it after time integration, so that altogether:
\begin{align}\label{BC+}
  (\sigma\vec{n})_{\rm t}+\kappa(\pi)\vv_{\rm t}=0,\ \ \vv_{\rm n}=0,
  \ \text{ and }\
  \frac1{H}\nabla\pi\cdot\vec{n}+\frac\pi c=0
\end{align}
with $c>0$ an elastic modulus related with an response of the boundary wall,
and
$\pi_{\rm ext}^{}$ a prescribed external pressure. This gives the boundary
contribution $\int_\Gamma\frac1{2c}\pi^2\,\d S$ into the stored energy.
More in detail, the calculus \eq{calculus-pi} enhances as
\begin{align}\nonumber
  \int_\Omega\nabla\pi\cdot \vv\,\d x&
  -\int_\Gamma\pi\,(\vv\cdot\vec{n})\,\d S
  \\[-.5em]&\nonumber=-\int_\Omega\pi\,{\rm div}\,\vv\,\d x=
    \int_\Omega\pi\Big(\frac1K\DT\pi-\frac1H\Delta\DT\pi\Big)\,\d x
    \\&\nonumber=\frac{\d}{\d t}\int_\Omega\frac1{2K}|\pi|^2
    +\frac1{2H}|\nabla\pi|^2\,\d x
    -\int_\Gamma\frac1{H}\pi(\nabla\DT\pi\cdot\vec{n})\,\d S
    \\&=\frac{\d}{\d t}\bigg(\int_\Omega\frac1{2K}|\pi|^2
    +\frac1{2H}|\nabla\pi|^2\,\d x
    +\int_\Gamma
   \frac{\pi^2}{2c} \,\d S\bigg)\,.
  \label{calculus-pi+}\end{align}
The energetics \eq{energy} is now modified as
\begin{align}\nonumber
&\int_\Omega\!\!\!\!\linesunder{\frac\varrho2|\vv(t)|^2}{kinetic}{energy}\!\!\!\!\!+
  \!\!\!\!\!\linesunder{
    \frac1{2K}\pi(t)^2
    +\frac1{2H}|\nabla\pi(t)|^2}{bulk stored}{energy}\!\!\!\!\!\,\d x
  +\int_\Gamma\!\!\!\!\!\!\!\!\!\!\linesunder{\frac{\pi(t)^2}{2c}}{boundary wall}{stored energy}\!\!\!\!\!\!\!\!\!\d S
\\&\nonumber\qquad
  +\int_0^t\!\!\int_\Omega\!\!\!\!\!\!\linesunder{\bbD e(\vv){:}e(\vv)
  }{dissipation}{by viscosity}
  \!\!\!\!\d x\d t+\int_0^t\!\!\int_\Gamma\!\!\!\!\!\!\!\!\!\!\!\!\!\!\!\!\!\linesunder{\kappa|v_{\rm t}|^2}{dissipation by a}{boundary ``friction''}\!\!\!\!\!\!\!\!\!\!\!\!\!\!\!\!\!\,\d S\d t=\int_0^t\!\!\int_\Omega\!\!\!\!\!\!\!\linesunder{g\cdot \vv}{power}{of loading}\!\!\!\!\!\!\!\d x\d t
  \\&\qquad\qquad\qquad\qquad
  +\int_\Omega\frac\varrho2|\vv_0|^2+\frac1{2K}\pi_0^2+\frac1{2H}|\nabla\pi_0|^2
  \,\d x+\int_\Gamma
  \frac{\pi_0^2}{2c}\,\d S
  \,.
  \label{energy-gradient}\end{align}

\begin{proposition}[Existence of weak solutions to \eq{NS-a}
    with \eq{Sobolev}]\label{prop-1}
  Let $d\le3$ and let
\begin{align}\label{ass}
 \!\!g\!\in\! L^1(I;L^2(\Omega;\R^d))+L^2(I;L^{6/5}(\Omega;\R^d)),\ \ 
  \vv_0\!\in\!L^2(\Omega;\R^d),\,\text{ and }\,\pi_0\!\in\!H^1(\Omega).
\end{align}
  Then:\\
  \Item{(i)}{the system \eq{NS-a} and \eq{Sobolev} with the initial and
  boundary conditions \eq{IC} and \eq{BC+}
  has a weak solution
  $\vv\in C_{\rm w}(I;L^2(\Omega;\R^d))\:\cap\:L^2(I;H^1(\Omega;\R^d))$
  and $\pi\in H^1(I;H^1(\Omega))$.}
\Item{(ii)}{If $\Omega$ is smooth, then even $\pi\in H^1(I;H^2(\Omega))$ so
  that,
in particular $\pi\in C(I{\times}\barOmega)$.}
\Item{(iii)}{If $d=2$, then the weak solution is unique.}
\end{proposition}

  \noindent{\it Sketch of the proof.}
  We can use some approximation, say by the Faedo-Galerkin method, using
  a nested sequences of finite-dimensional subspaces of
  $H^1(\Omega;\R^d)$ and of $H^1(\Omega)$ whose union is dense in
  these Sobolev spaces to be used for the approximation of
  $\vv$ and of $\pi$, respectively. Let us index them by $k\in\N$ and denote
  by $(\vv_k,\pi_k)$ the approximate solutions created by this way.
  In particular, the discrete variant of 
  \eq{Sobolev} leads to the integral identity
  \begin{align}\label{Sobolev-k}
    \int_\Omega\frac1K\DT\pi_k\widetilde\pi+
    \Big(\frac1H\nabla\DT\pi_k-\vv_k\Big)\cdot\nabla\widetilde\pi\,\d x=0
  \end{align}
  for all $\widetilde\pi$ from the mentioned finite-dimensional subspace
  of $H^1(\Omega)$ so that $\widetilde\pi=\pi_k$ will be a legitimate
  test to be used for the calculus \eq{calculus-pi+} which then
  holds for $(\vv_k,\pi_k)$ too.
  
  The local-in-time existence of $(\vv_k,\pi_k)$ is due to the standard
  arguments from the theory of ordinary differential systems together with
  the subsequent prolongation
  on the whole time interval $I$ due to the $L^\infty(I)$-estimates
  which follows from \eqref{energy-gradient} which holds also for
  $(\vv_k,\pi_k)$. These
  estimates are obtained from \eqref{energy-gradient} by Young's and Gronwall's
  inequalities and by using the embedding $H^1(\Omega;\R^d)\subset L^6(\Omega;\R^d)$
  which holds for $d\le3$. More specifically, we can see the a-priori estimates
  \begin{subequations}\label{est-gradient}\begin{align}\label{est1-gradient}
    &\|\vv_k\|_{L^\infty(I;L^2(\Omega;\R^d))\:\cap\:L^2(I;H^1(\Omega;\R^d))}^{}\le C\,,
    \\&\label{est2-gradient}\|\pi_k\|_{L^\infty(I;H^1(\Omega))}^{}\le C\,.
\end{align}\end{subequations}
  By the Banach selection principle, we can select a  subsequence
  converging weakly* to a limit $(u,\pi)$. This limit is a weak solution
  to the mentioned initial-boundary-value problem for \eq{NS-a} and \eq{Sobolev}.
  The limit passage in the Galerkin approximation for $k\to\infty$ is simply due to the
  mentioned weak* convergence. For the only two
  nonlinear terms in \eq{NS-a} one can use the strong convergence of $\vv_k$
   in $L^{10/3-\epsilon}(\Omega;\R^d)$ with $\epsilon>0$ arbitrarily
  small. This follows (through interpolation) 
by  the Aubin-Lions compact-embedding theorem,
using also a uniform bound of a Hahn-Banach extension of
$\varrho\DT\vv_k$  on $L^1(I;B)$ with an arbitrary Banach space $B$,
here  on the space from the former inclusion in \eqref{est+},
  cf.\ \cite[Sect.8.4]{Roub13NPDE}
  for details about this technique.
  
  Further, testing \eq{Sobolev-k}
  by $\DT\pi_k$ yields $\DT\pi_k$ bounded in $L^2(I;H^1(\Omega))$. This is
  inherited in the limit $\DT\pi\in L^2(I;H^1(\Omega))$.

For a smooth $\Omega$, we can use $H^2$-regularity of
the $\Delta$-operator so that, from
$\Delta\DT\pi=H(\DT\pi/K+{\rm div}\,\vv)\in L^2(I;L^2(\Omega))$,
we obtain also $\DT\pi\in L^2(I;H^2(\Omega))$.


As for the uniqueness in the 2-dimensional case,
we just expand the estimation which, in the incompressible case,
are based on the Gagliardo-Nirenberg
  inequality $\|\vv\|_{L^4(\Omega;\R^d)}\le C\|\vv\|_{L^2(\Omega;\R^d)}^{1/2}\|\nabla\vv\|_{L^2(\Omega;\R^{d\times d})}^{1/2}$.
Abbreviating $\vv_{12}:=v_1{-}v_2$ and $\pi_{12}:=\pi_1{-}\pi_2$ for two weak
solutions $(\vv_1,\pi_1)$ and $(\vv_2,\pi_2)$, one has the estimate
\begin{align}\nonumber&
\frac{\d}{\d t}\Big(\frac\varrho2\|\vv_{12}\|_{L^2(\Omega;\R^d)}^2
+\frac1{2K}\|\pi_{12}\|_{L^2(\Omega)}^2
+\frac1{2H}\|\nabla\pi_{12}\|_{L^2(\Omega;\R^d)}^2\Big)
\\\nonumber&\hspace{21em}
+\min(K_\text{\sc v}^{},G_\text{\sc v}^{})\|e(\vv_{12})\|_{L^2(\Omega;\R^{d\times d})}^2
\\\nonumber&\hspace{0em}
=\int_\Omega\varrho\Big((\vv_1{\cdot}\nabla)\vv_1+\frac12({\rm div}\,\vv_1)\,\vv_1
-(\vv_2{\cdot}\nabla)\vv_2-\frac12({\rm div}\,\vv_2)\,\vv_2\Big){\cdot}\vv_{12}
\,\d x\!
\\\nonumber&\hspace{0em}
=\int_\Omega\varrho\Big((\vv_{12}{\cdot}\nabla)\vv_1+\frac12({\rm div}\,\vv_{12})\,\vv_2
+(\vv_2{\cdot}\nabla)\vv_{12}+\frac12({\rm div}\,\vv_1)\,\vv_{12}\Big)\cdot \vv_{12}
\,\d x
\\\nonumber
&\hspace{0em}\le
\!\frac{3\varrho}2\|\nabla\vv_1\|_{L^2(\Omega;\R^{d\times d})}^{}\|\vv_{12}\|_{L^4(\Omega;\R^d)}^2
\!+
\!\frac{3\varrho}2
\|\vv_2\|_{L^4(\Omega;\R^d)}^{}\|\nabla\vv_{12}\|_{L^2(\Omega;\R^{d\times d})}^{}\|\vv_{12}\|_{L^4(\Omega;\R^d)}^{}\
\\&\nonumber
\le
C
\|\nabla\vv_1\|_{L^2(\Omega;\R^{d\times d})}^{}\|\vv_{12}\|_{L^2(\Omega;\R^d)}^{}\|\nabla\vv_{12}\|_{L^2(\Omega;\R^d)}^{}
\\&\nonumber\hspace{15em}
+C\|\vv_2\|_{L^4(\Omega;\R^d)}^{}\|\nabla\vv_{12}\|_{L^2(\Omega;\R^{d\times d})}^{3/2}\|\vv_{12}\|_{L^2(\Omega;\R^d)}^{1/2}\
\\&\nonumber\le
\epsilon\|\nabla\vv_{12}\|_{L^2(\Omega;\R^{d\times d})}^{2}
+C^2\|\nabla\vv_1\|_{L^2(\Omega;\R^{d\times d})}^2\|\vv_{12}\|_{L^2(\Omega;\R^d)}^2/\epsilon
\\&\hspace{15em}
+
C^4\|\vv_2\|_{L^4(\Omega;\R^d)}^4\|\vv_{12}\|_{L^2(\Omega;\R^d)}^{2}/\epsilon^3\,,
\label{NS-uniqueness+}\end{align}
where the constant $C$ is from  the Gagliardo-Nirenberg inequality.
Then taking $\epsilon>0$ sufficiently small and using Korn's and
Gronwall's inequalities,
we can see that $\vv_{12}=0$ and $\pi_{12}=0$. For this, it is also
important that $\vv_2\in L^4(I{\times}\Omega;\R^d)$, which follows from
the interpolation 
$L^\infty(I;L^2(\Omega))\cap L^2(I;H^1(\Omega))\subset L^4(I;L^4(\Omega))$
holding for $d=2$ again by the Gagliardo-Nirenberg inequality.
\hfill$\Box$

\medskip

\def\Kv{D}

In comparison with Sect.~\ref{sect-simple}, the nondissipative-gradient model brings
an additional mechanism contributing to {\it dispersion} not related with any
{\it attenuation}. To see a more specific character of the
dispersion-attenuation behind \eqref{Sobolev}, let us neglect the nonlinear
convective terms in \eqref{NS-a} and consider a 1-dimensional
 linearized  situation. We thus consider the simple system
\begin{align}\label{1D-system}
  &\DT u=\vv,\ \ \ \ \ \ \ \ \varrho\DT\vv-\Kv v_{xx}^{}+\pi_x^{}=0,\ \ \ \text{ and }\ \ \ 
  \frac{\DT\pi}K+\vv_x^{}=\frac{\DT\pi_{xx}^{}}H\,.
\end{align}
Equivalently, we can also write it as $\varrho\DDT u-\Kv\DT u_{xx}+\pi_x^{}=0$
and $\frac1K\pi+u_x^{}=\frac1H\pi_{xx}^{}$. In particular,
we eliminate $\vv$ and consider only the displacement $u$ and the pressure to have 
\begin{align*}
  \frac\varrho K\DDT u-\frac{\Kv}K\DT u_{xx}+\frac{\pi_x^{}}K=0\ \ \ \text{ and }\ \ \
   \frac\varrho H\DDT u_{xx}^{}-\frac{\Kv}H\DT u_{xxxx}^{}+\frac{\pi_{xxx}^{}}H=0\,.
\end{align*}
When subtracting these two equations, we can eliminate also $\pi$ by using
$\frac1K\pi_x^{}-\frac1H\pi_{xxx}^{}=-u_{xx}^{}$
and obtain a dispersive equation
\begin{align}\label{dispersion}
  \varrho \DDT u-\varepsilon\varrho\DDT u_{xx}^{}
  -\Kv \DT u_{xx}^{}+\varepsilon \Kv\DT u_{xxxx}^{}
  -Ku_{xx}^{}=0\,\ \ \ \ \text{ with }\ \varepsilon:=\frac KH
  \,.
\end{align}
The term $\varepsilon\varrho\DDT u_{xx}^{}$ arisen in \eqref{dispersion} is
sometimes called {\it micro-inertia}; for various models involving such
micro-inertia and its role in dispersion we refer in particular to
\cite[Chap.~6]{BerVan17IVT} and
\cite{AskAif09GEFW,EBPB05WMMD,LaMaAi06TNEB,MNAB17RMIE}.
Analogously, $\varepsilon \Kv\DT u_{xxxx}$ might be called micro-viscosity.
A conventional way to calculate dispersion and attenuation in linear
media is to use the ansatz $u={\rm e}^{{\rm i}(w t+x/\lambda)}$ with the angular
frequency $w=\omega+{\rm i}\gamma$ considered complex with 
$\omega,\gamma\in\R$ and the real-valued 
wavelength $\lambda$; here ${\rm i}=\sqrt{-1}$ denotes the imaginary unit.
Then $\varrho\DDT u=-\varrho w^2u$,
$\varepsilon\varrho\DDT u_{xx}^{}=\varepsilon\varrho w^2u/\lambda^2$,
$\Kv\DT u_{xx}=-{\rm i}\Kv w u/\lambda^2$,
$\varepsilon \Kv\DT u_{xxxx}={\rm i}\varepsilon \Kv w u/\lambda^4$, and 
$Ku_{xx}=-Ku/\lambda^2$, so that the one-dimensional dispersive wave equation
\eq{dispersion} yields the algebraic condition 
\begin{align}\label{KK-relation-}
\varrho\eta
w^2
-{\rm i}\Kv
\eta\frac w{\lambda^2}
-\frac K{\lambda^2}=0\ \ \ \text{ with }\
\eta:=1+\frac{\varepsilon}{\lambda^2}\,.
\end{align}
When substituting $w=\omega+{\rm i}\gamma$ so that $w^2=\omega^2-\gamma^2
+2{\rm i}\omega\gamma$, we obtain two algebraic equations for the 
real and the imaginary part each, called {\it Kramers-Kronig's 
relations}. More specifically, here
\begin{align}\label{KK-relation}
&\varrho\eta(\omega^2-\gamma^2)=\frac K{\lambda^2}-\Kv\eta\frac\gamma{\lambda^2}
\ \ \ \text{ and }\ \ \
2\varrho\gamma=\frac{\Kv}{\lambda^2}\,.
\end{align}
Rather exceptionally, these algebraic relations can explicitly be solved: From
the latter equation we can read that $\gamma=\Kv/(2\varrho\lambda^2)$
and then the former equation yields $\omega^2=K/(\varrho\eta\lambda^2)
-\Kv^2/(4\varrho^2\lambda^4)$. When substituting $\eta$ from \eq{KK-relation-},
we obtain the ({\it real} part of the) {\it angular frequency}
\begin{align}\label{omega=omega(lambda)}
  \omega=\omega(\lambda)=\sqrt{\frac K{\varrho\varepsilon{+}\varrho\lambda^2}
  -\frac{\Kv^2}{4\varrho^2\lambda^4}}\,.
\end{align}
Then, realizing that the speed of
waves is $v=\omega\lambda$ and substituting $\eta$ from \eq{KK-relation-}, we
obtain 
\begin{align}\label{v=v(lambda)}
\vv=\vv(\lambda)=\sqrt{\frac K{\varrho(1+\varepsilon/\lambda^2)}-\frac{\Kv^2}{4\varrho^2\lambda^2}}\,,
\end{align}
which gives a {\it normal dispersion}\index{dispersion!normal} for sufficiently
long waves, namely having a length $\lambda>\lambda_\text{\sc crit}$ with the
critical wave length solving the equation
\begin{align}\label{crit}
  4\varrho K\lambda_\text{\sc crit}^4
  -\Kv^2\lambda_\text{\sc crit}^2-\varepsilon \Kv^2=0\,.
\end{align}
In terms of the so-called {\it angular wavenumber} $k=2\uppi/\lambda$,
\eq{v=v(lambda)} can be written as
\begin{align}\label{v=v(k)}
\vv=\vv(k)=\sqrt{\frac K{\varrho(1+\varepsilon k^2/(4\uppi^2))}-\frac{\Kv^2k^2}{16\varrho^2\uppi^2}}\,,
\end{align}
cf.\ Fig.~1. Let us recall that 
the adjective ``normal'' for dispersion means that waves with longer lengths 
propagate faster than those with shorter lengths. Waves with the 
length $\lambda_\text{\sc crit}$ or shorter are so fast attenuated that 
they cannot propagate at all. Let us note that the coefficient $\gamma$ 
determining the attenuation is here inversely proportional to the square power 
of the wave length.
This also reveals the dispersion/attenuation for the simple model from
Sect.~\ref{sect-simple} when putting $\varepsilon=0$ into \eq{v=v(lambda)}
and \eq{crit}; in particular $\lambda_\text{\sc crit}=\Kv/\sqrt{4\varrho K}$.
The other extreme case is for the inviscid situation $\Kv\to0+$ where
$\lambda_\text{\sc crit}\to0+$ so that also waves with ultra-high frequencies
can propagate. Let us note that the coefficient $\gamma$ determining the
attenuation is here inversely proportional to the square power of the wave
length, i.e.\ naturally the attenuation rises for high frequencies.

From \eq{omega=omega(lambda)}, we can also read the so-called 
{\it group velocity} of waves, which is the velocity with which the overall
envelope shape of the wave's amplitudes propagates through space and
which is defined as ${\rm d}\omega/{\rm d}k$ with $k=2\uppi/\lambda$.
Here, \eq{omega=omega(lambda)}
with $\lambda=2\uppi/k$ yields
\begin{align*}
  \omega=\omega(k)=\sqrt{\frac{K k^2}{4\uppi^2\varrho{+}\varepsilon\varrho k^2}
  -\frac{\Kv^2k^4}{64\uppi^2\varrho^2}}\,.
\end{align*}
{}In this case, 
we obtain the group velocity:
\begin{align}\nonumber
  \frac{{\rm d}\omega}{{\rm d}k}=\frac1{\omega(k)}\bigg(\frac{Kk}{\sqrt{4\uppi^2\varrho{+}\varepsilon\varrho k^2}}
    -\frac{K\varepsilon\varrho k^3}{\sqrt[3/2]{4\uppi^2\varrho{+}\varepsilon\varrho k^2}}
-\frac{D^2k^3}{32\uppi^2\varrho^2}\bigg)\,.
\end{align}

\begin{center}\hspace*{-2em}
\begin{tikzpicture}[scale=0.55]
\begin{axis}[
    scale only axis,
    grid=major,
    axis lines=middle,
    inner axis line style={=>},
    xlabel={$\text{\Huge $k$}$\hspace*{-40em}},
    ylabel={\Huge $\hspace*{1em}v{=}v(k)$},
    ytick={0,1,...,5},
    xtick={0,1,...,5},
    ymin=0,
    ymax=5,
    xmin=-.0,
    xmax=5,
]
    \pgfplotsset{samples=100} 

\addplot +[mark=none,black,very thick] {sqrt(max(0,9/(1+.1*x^2/39.5)-9*x^2/158))};
\addplot +[mark=none,black,very thick,densely dashed] {sqrt(max(0,9/(1+1*x^2/39.5)-9*x^2/158))};
\addplot +[mark=none,black,very thick,loosely dashed] {sqrt(max(0,9/(1+10*x^2/39.5)-9*x^2/158))};
 
    \legend{$v$ from \eqref{v=v(k)} for $\varepsilon=0.1$,
      $v$ from \eqref{v=v(k)} for $\varepsilon=1$,
      $v$ from \eqref{v=v(k)} for $\varepsilon=10$,
            };
\end{axis}
\end{tikzpicture}\hspace*{-1.5em}
\begin{tikzpicture}[scale=0.55]
\begin{axis}[
    scale only axis,
    grid=major,
    axis lines=middle,
    inner axis line style={=>},
    xlabel={$\text{\Huge $\lambda$}$\hspace*{-37em}},
    ylabel={$\text{\Huge $\hspace*{.8em}v{=}v(\lambda)\hspace*{-5.4em}$}$},
     ytick={0,1,...,5},
    xtick={0,1,...,5},
    ymin=0,
    ymax=5,
    xmin=-.0,
    xmax=5,
]
    \pgfplotsset{samples=100} 
    \addplot +[mark=none,black,very thick] {sqrt(max(0,9/(1+.1/x^2)-9/(4*x^2)))};
    \addplot +[mark=none,densely dashed,black,very thick] {sqrt(max(0,9/(1+1/x^2)-9/(4*x^2)))};
    \addplot +[mark=none,loosely dashed,black,very thick] {sqrt(max(0,9/(1+10/x^2)-9/(4*x^2)))};
   
    \legend{$v$ from \eqref{v=v(lambda)} for $\varepsilon=0.1$,
      $v$ from \eqref{v=v(lambda)} for $\varepsilon=1$,
      $v$ from \eqref{v=v(lambda)} for $\varepsilon=10$,
            };
\end{axis}
\end{tikzpicture}\hspace*{1em}
   
\nopagebreak
{\small\it
\hspace*{-2em}Fig.\,1:~\begin{minipage}[t]{37em}
Dependence of the velocity of sinusoidal waves on
  the angular wavenumber $k$ (left) and 
  on the wave length $\lambda=2\uppi/k$ (right);
  an illustration of the normal dispersion due to \eq{v=v(k)}
  and \eq{v=v(lambda)} for $K$=9 and $\varrho=1$, and $\Kv=3$.
  Waves with ultra short lengths (or with ultra high wave numbers) have zero velocity, i.e.\ cannot propagate.
\end{minipage}}\end{center}

There is an alternative way how to incorporate a conservative gradient
term into the model: instead of enhancing \eq{NS-b} to \eq{Sobolev},
we can enhance the pressure term in \eq{NS-a}. More specifically,
the system \eq{NS} is enhanced for 
\begin{subequations}\label{NS-Delta-pi}
  \begin{align}
    &\varrho\DT \vv+\varrho(\vv{\cdot}\nabla)\vv-{\rm div}
\big(\bbD e(\vv)-(\pi{-}\ell^2\Delta\pi)\bbI\big)
=g-\frac\varrho2({\rm div}\,\vv)\,\vv
\label{NS-Delta-pi-a}
\,,
    \\\label{NS-Delta-pi-b}
    &\DT\pi+K\,{\rm div}\,\vv\,=0\,
\end{align}\end{subequations}
with some $\ell>0$ presumably small, having the physical
dimension in meters and determining a certain length scale for possible
spacial pressure variations. Modifying still the boundary conditions \eq{BC+} as
\begin{align}\label{BC++}
  (\sigma\vec{n})_{\rm t}+\kappa(\pi)\vv_{\rm t}=0\,,\ \ \ \vv_{\rm n}=0\,,
  \ \text{ and }\
  \frac1{K}\nabla\pi\cdot\vec{n}+\frac{\DT\pi}c=0\,,
\end{align}
with $\sigma=\bbD e(\vv)-(\pi{-}\ell^2\Delta\pi)\bbI$,
the previous calculus \eq{calculus-pi+} now modifies as
\begin{align}\nonumber
  \int_\Omega\nabla(\pi{-}\ell^2\Delta\pi)\cdot \vv\,\d x&
  -\int_\Gamma(\pi{-}\ell^2\Delta\pi)\,(\vv\cdot\vec{n})\,\d S
\\&\nonumber =
  -\int_\Omega(\pi{-}\ell^2\Delta\pi)\,{\rm div}\,\vv\,\d x=
    \int_\Omega(\pi{-}\ell^2\Delta\pi)\frac1K\DT\pi\,\d x
    \\&\nonumber=\frac{\d}{\d t}\int_\Omega\frac1{2K}\pi^2
    +\frac{\ell^2}{2K}|\nabla\pi|^2\,\d x
    -\int_\Gamma\frac{\ell^2}K\pi(\nabla\pi\cdot\vec{n})\,\d S
   \\&=\frac{\d}{\d t}\bigg(\int_\Omega\frac1{2K}\pi^2
    +\frac{\ell^2}{2K}|\nabla\pi|^2\,\d x+
    \int_\Gamma\frac{\ell^2}{2c}\pi^2\,\d S\bigg)\,.
\label{calculus-pi++}\end{align}
This gives again the energy balance \eq{energy-gradient} only with
$H$ replaced by $\ell^2/K$ and with $c$ replaced by $c/\ell^2$.
The Proposition~\ref{prop-1}(i) and (iii) holds for
the boundary-value problem \eq{NS-Delta-pi}--\eq{BC++} with the initial
condition \eq{IC}, too.

Although the model \eq{NS-Delta-pi}--\eq{BC++} is energetically and
analytically very similar to \eq{NS-a} with \eq{Sobolev},
it differs as far as the dispersive character. Indeed, 
instead of \eq{1D-system} we now have
\begin{align}\label{1D-system-Delta-pi}
  &\DT u=\vv,\ \ \ \ \ \ \ \ \varrho\DT\vv-\Kv v_{xx}^{}+\pi_x^{}
  -\ell^2\pi_{xxx}^{}=0,\ \ \ \text{ and }\ \ \ \DT\pi+K\vv_x^{}=0\,.
\end{align}
The last equation means $\pi+Ku_x^{}=0$, so that \eq{1D-system-Delta-pi}
results, instead of \eq{dispersion}, to the dispersive equation:
\begin{align}\label{disp-eq}
  \varrho \DDT u-\Kv\DT u_{xx}^{}-K u_{xx}^{}+\ell^2 K u_{xxxx}^{}=0\,.
\end{align}
Instead of \eq{KK-relation-}, we now have $\varrho w^2
-{\rm i}\Kv w/\lambda^2-K/\lambda^2-\ell^2 K/\lambda^4=0$
so that the Kramers-Kronig's relations look as
\begin{align}\label{KK-relation+}
&\varrho(\omega^2-\gamma^2)=\frac K{\lambda^2}-\Kv\frac\gamma{\lambda^2}+\ell^2 \frac K{\lambda^4}
\ \ \ \text{ and }\ \ \
2\varrho\gamma=\frac{\Kv}{\lambda^2}\,.
\end{align}
The speed of wave $\vv=\omega\lambda$ is now
\begin{align}\label{v-anomalous}
v=v(\lambda)=\sqrt{\frac K\varrho+\frac{\ell^2 K}{\varrho\lambda^2}-\frac{\Kv}{4\varrho^2\lambda^2}}\,.
\end{align}
In terms of the angular wavenumber $k=2\uppi/\lambda$, we now have,
\begin{align}\label{v-anomalous+}
  v=v(k)=\sqrt{\frac K\varrho+\Big(\frac{\ell^2 K}{4\uppi^2\varrho}
    -\frac{\Kv}{16\uppi^2\varrho^2}\Big)k^2}\,.
\end{align}
If the dissipative effects do not dominate (i.e.\ $\Kv>0$ is
small), \eq{v-anomalous} gives an {\it anomalous dispersion}, i.e.\
higher-frequency waves (i.e.\ with longer wave length) propagate faster
than waves with lower frequencies, cf.\ Fig.~2. From the term $\ell^2 K u_{xxxx}^{}$ in
\eq{disp-eq}, we can see that the model \eq{NS-Delta-pi} actually
deal with the strain gradient while the previous model
\eq{NS-a} with \eq{Sobolev} has a gradient of stress (here pressure).
This is respectively the Aifantis' versus Eringen's approach to incorporation
of internal length scale, cf.\ \cite{Aifa92RGLD,AskAif09GEFW,AskAif11GESD,Erin83DENE}.
\begin{center}\hspace*{-2em}
\begin{tikzpicture}[scale=0.55]
\begin{axis}[
    scale only axis,
    grid=major,
    axis lines=middle,
    inner axis line style={=>},
    xlabel={$\text{\Huge $k$}$\hspace*{-40em}},
    ylabel={\Huge $\hspace*{1.em}v{=}v(k)$},
    ytick={.5,1,...,2},
    xtick={0,1,...,5},
    ymin=.5,
    ymax=2,
    xmin=-.0,
    xmax=5,
]
    \pgfplotsset{samples=100} 
    \addplot +[mark=none,black,very thick] {1*sqrt(1+x^2/39.5)};
    \addplot +[mark=none,densely dashed,black,very thick] {1*sqrt(1+x^2/247)};
    \addplot +[mark=none,loosely dashed,black,very thick] {1*sqrt(1+x^2/3950)};
   
    \legend{$v$ from \eqref{v-anomalous+} for $\ell=1$,
      $v$ from \eqref{v-anomalous+} for $\ell=.4$,
      $v$ from \eqref{v-anomalous+} for $\ell=.1$,
            };
\end{axis}
\end{tikzpicture}\hspace*{-1.5em}
\begin{tikzpicture}[scale=0.55]
\begin{axis}[
    scale only axis,
    grid=major,
    axis lines=middle,
    inner axis line style={=>},
    xlabel={$\text{\Huge $\lambda$}$\hspace*{-37em}},
    ylabel={$\text{\Huge $\hspace*{.8em}v{=}v(\lambda)\hspace*{-5.4em}$}$},
     ytick={0,1,...,3},
    xtick={0,1,...,5},
    ymin=0,
    ymax=3,
    xmin=-.0,
    xmax=5,
]
    \pgfplotsset{samples=100} 
       \addplot +[mark=none,black,very thick] {1*sqrt(1+1/x^2)};
    \addplot +[mark=none,densely dashed,black,very thick] {1*sqrt(1+.16/x^2)};
    \addplot +[mark=none,loosely dashed,black,very thick] {1*sqrt(1+.01/x^2)};
  
    \legend{
      $v$ from \eqref{v-anomalous} for $\ell=1$,
      $\,v$ from \eqref{v-anomalous} for $\ell=.4$,
      $\ v$ from \eqref{v-anomalous} for $\ell=.1$,
             };
\end{axis}
\end{tikzpicture}\hspace*{1em}

\nopagebreak
{\small\it
\hspace*{-2em}Fig.\,2:~\begin{minipage}[t]{37em}
Dependence of the velocity of sinusoidal waves on
  the angular wavenumber $k$ (left) and 
  on the wave length $\lambda=2\uppi/k$ (right);
  an illustration of the anormalous dispersion due to \eq{v-anomalous+}
  and \eq{v-anomalous} for $K$=1, $\varrho=1$, and $\Kv>0$ very small.
\end{minipage}}\end{center}

\begin{remark}[{\sl Pressure-dependent compressibility}]\label{rem-K=K(pi)}
  \upshape
  Actually, the compressibility and thus also the speed of P-waves may also depend
on the pressure, cf.\ e.g.\ \cite[Table~V]{FinMil73CWFT} or \cite{MCBS80NHPE}
for the case of seawater. It is not difficult to make $K$ dependent on $\pi$
in a general continuous manner by replacing the term $\frac1{2K}|\pi|^2$
by a general function $\phi(\pi)$, which leads to $K(\pi)=\pi/\phi'(\pi)$.
Indeed, the calculus \eq{calculus-pi+} thus modifies by
using $\pi\DT\pi/K(\pi)=\phi'(\pi)\DT\pi=\frac{\partial}{\partial t}\phi(\pi)$.
\end{remark}

\begin{remark}[{\sl A generalization of the boundary conditions}]
\upshape  
  As the pressure is well defined on the boundary, one can think about a
  penetrable wall with a flux depending of pressure, which could be modeled
  by the boundary condition $\vv_{\rm n}=\vv\cdot\vec{n}=f(\pi)$. Yet, it would
  expand the energy balance \eq{energy}  by the boundary flux
  $\int_0^t\!\!\int_\Gamma\frac\varrho2|\vv|^2f(\pi)\d S\d t$ which
  seems to corrupt the a-priori estimation even if $f(\pi)$ would be
  supposed bounded.
\end{remark}

\section{Towards dispersion through dissipative gradient terms}\label{sec-dispersive2}

An interesting modification has been devised by A.P.\,Oskolkov
\cite{Osko73SPQL} who considered (rather for analytical purposes) a parabolic
regularization of the incompressibility constraint of the type:
\begin{align}\label{Oskolkov}
\frac1K\DT\pi+{\rm div}\,\vv=\frac1H\Delta\pi\,,
\end{align}
where the bulk ``hyper-modulus'' $H$ has now the physical dimension
Pa\,s/m$^2$. Ignoring the convective term, this was later developed under the
name ``mixed quasi-compressibility methods'' in
\cite[Ch.\,5]{Proh97PQCM}. It actually combines the
simple time-regularization \eq{NS-b} with a mere elliptic
regularization ${\rm div}\,\vv=\frac1H\Delta\pi$ as used e.g.\ in
\cite{BerSpi17CSWS,BuMaRa07NSEN,Proh97PQCM}.

It can be derived from a fully compressible model when introducing
some {\it diffusion} into the {\it convective mass transport} \eq{compress-NS1},
i.e.\ 
\begin{align}\label{Oskolkov-}
  \DT\rho+
  \vv{\cdot}\nabla\rho+\rho\,{\rm div}\,\vv
  =\epsilon\Delta\rho\qquad\text{ with }\ \ \epsilon=K/H
\end{align}
for some hyper-modulus $H$. Although it does not seem
 deducible  from the
Boltzmann equation by usual procedure toward gas-dynamic limit
\cite{BaGoLe91FDLK}, such a convection-diffusion transport is occasionally
considered, devised by H.\,Brenner \cite{Bren05KVT,Bren06FMR},
cf.\ also \cite{BarOtt07CBMN,FeiVas09NPFD}, although sometimes rather as an
artificial
regularization only, cf.\ \cite[Eq.\,(3.172)]{FeiNov09SLTV}. Similarly,
\cite{BuMaMi17SDRN} devised a diffusion of the deviatoric part of stress
in the incompressible case as a regularization.
This falls into a general concept of a {\it parabolic perturbation} of the
(here mass-) conservation law. For a thermodynamical justification
see also \cite{VaPaGr17EMFF}, opposing \cite{OtStLi09IDCM}.
Then the small-perturbation ansatz \eq{ansatz} would lead, instead of
\eq{NS-b}, to 
\eq{Oskolkov}.
The physical dimension of $\epsilon$ is m$^2$/s. Vaguely speaking,
dividing $\epsilon$ by a ``characteristic'' velocity of the flow and
a ``characteristic'' length of the system, we obtain a dimensionalless 
{\it P\'eclet number} expressing dominance of either the convective or the
diffusive transport phenomena.

Diffusion processes are always dissipative, so it is not surprising that this
``regularization'' of the incompressibility condition is dissipative, too.
The calculus \eq{calculus-pi} now enhances as
\begin{align}\nonumber
  \int_\Omega\nabla\pi\cdot \vv\,\d x&-\int_\Gamma\pi\,(\vv\cdot\vec{n})\,\d S=
  -\int_\Omega\pi\,{\rm div}\,\vv\,\d x=
    \int_\Omega\pi\Big(\frac1K\DT\pi-\frac1H\Delta\pi\Big)\,\d x
    \\[-.3em]&\ \ \ \ \ =\frac{\d}{\d t}\int_\Omega\frac1{2K}|\pi|^2\,\d x
    +\int_\Omega\frac1H|\nabla\pi|^2\,\d x
    -\int_\Gamma\frac1H\pi(\nabla\pi\cdot\vec{n})\,\d S\,.
  \label{calculus-pi+++}\end{align}
Considering again the boundary conditions \eq{BC+}, 
the energetics \eq{energy-gradient} now modifies as 
\begin{align}\nonumber
&\int_\Omega\!\!\!\!\linesunder{\frac\varrho2|\vv(t)|^2}{kinetic}{energy}\!\!\!\!\!+
  \!\!\!\!\!\linesunder{\frac1{2K}\pi(t)^2}{bulk stored}{energy}\!\!\!\!\!\,\d x
+\int_0^t\!\!\int_\Omega\!\!\!\!\!\!\linesunder{\bbD e(\vv):e(\vv)+\frac1H|\nabla\pi|^2
  }{dissipation by ``viscosity''}{in shear rate and pressure}
  \!\!\!\d x\d t
  \\&  \quad+\int_0^t\!\!\int_\Gamma\!\!\!\!\!\linesunder{
    \kappa|\vv_{\rm t}|^2+\frac1c\pi^2}{dissipation on}{the boundary}
  \!\!\!\!\!\d S\d t 
  =\int_\Omega\frac\varrho2|\vv_0|^2+\frac1{2K}\pi_0^2
  \,\d x+\int_0^t\!\!\int_\Omega\!\!\!\!\!\!\linesunder{g\cdot \vv}{power}{of loading}\!\!\!\!\!\!\!\d x\d t
  \,.
\label{energy+}\end{align}

\begin{proposition}[Existence of weak solutions to
   \eq{NS-a} with \eq{Oskolkov}]\label{prop-2}
  Let  $d\le3$ and the assumptions \eq{ass} hold.
  Then:\\
  \Item{(i)}{the system \eq{NS-a} and \eq{Oskolkov} with the initial and
  boundary conditions \eq{IC} and \eq{BC+}
  has a weak solution
  $\vv\in C_{\rm w}(I;L^2(\Omega;\R^d))\:\cap\:L^2(I;H^1(\Omega;\R^d))$
  and $\pi\in C_{\rm w}(I;H^1(\Omega))\:\cap\:H^1(I;L^2(\Omega))$.}
\Item{(ii)}{If $\Omega$ is smooth, then even $\pi\in L^2(I;H^2(\Omega))$ so
  that, in particular $\pi(t)\in C(\barOmega)$ for a.a.\ $t\in I$.}
\Item{(iii)}{If $d=2$, then the weak solution is unique.}
\end{proposition}

\noindent{\it  Sketch of the proof.}
  Again we use the Faedo-Galerkin method as in the proof of Proposition~\ref{prop-1}.
From \eq{energy+} written for the approximate solutions, we obtain a-priori
estimates 
\begin{subequations}\label{est-Oskolkov}\begin{align}\label{est1-Oskolkov}
&\|\vv_k\|_{L^\infty(I;L^2(\Omega;\R^d))\:\cap\:L^2(I;H^1(\Omega;\R^d))}^{}\le C\,,
\\&\label{est2-Oskolkov}\|\pi_k\|_{L^\infty(I;L^2(\Omega))\:\cap\:L^2(I;H^1(\Omega))}^{}
    \le C\,,
\intertext{In fact, \eq{est1-Oskolkov} is the former estimate
\eq{est1-gradient}. Further a-priori estimate follows by testing \eq{Oskolkov}
by $\DT\pi$, using that \eq{est1-Oskolkov} yields
${\rm div}\,\vv_k\in L^2(I{\times}\Omega)$. Thus we obtain}
&
  \|\pi_k\|_{L^\infty(I;H^1(\Omega))}{}\le C,\ \ \ \text{ and }\ \ \
  \label{est6-Oskolkov}
    \|\DT\pi_k\|_{L^2(I{\times}\Omega)}^{}\le C\,.
\end{align}\end{subequations}
By comparison, we obtain also an a-priori
information for a Hahn-Banach extension of $\varrho\DT\vv_k$
in $L^2(I;H^1(\Omega;\R^d)^*)+L^{5/4}(I{\times}\Omega;\R^d)$,
and then, by the Aubin-Lions theorem
used for such Hahn-Banach extension as in \cite[Sect.8.4]{Roub13NPDE},
we can pass to the limit as in the proof of Proposition~\ref{prop-1}.

As to (ii), assuming also $\Omega$ to
be smooth so that $H^2$-regularity of the Laplacian is at disposal,
from
$\Delta\pi=H({\rm div}\,\vv+\DT\pi/K)\in L^2(I{\times}\Omega)$ we obtain the
$L^2(I;H^2(\Omega;\R^d))$-estimate of $\pi$. 


As for (iii), we just modify the estimate \eq{NS-uniqueness+} appropriately:
 more specifically,
$\frac{\d}{\d t}\frac1{2H}\|\nabla\pi_{12}\|_{L^2(\Omega;\R^d)}^2$
is to be replace{d}  by $\frac1{H}\|\nabla\pi_{12}\|_{L^2(\Omega;\R^d)}^2$
otherwise the argumentation \eq{NS-uniqueness+} holds unchanged.
\hfill$\Box$

\medskip

The 1-dimensional calculations towards dispersion/attenuation of the P-waves
now modify \eq{1D-system} as
\begin{align}\label{1D-system+}
  &\DT u=\vv,\ \ \ \ \ \ \ \ \varrho\DT\vv-\Kv v_{xx}+\pi_x^{}=0,\ \ \ \text{ and }\ \ \ 
  \frac{\DT\pi}K+\vv_x^{}=\frac{\pi_{xx}^{}}H\,.
\end{align}
Eliminating $v$, we can write it as $\varrho\DDT u-\Kv\DT u_{xx}+\pi_x^{}=0$
and $\frac1K\DT\pi+\DT u_x^{}=\frac1H\pi_{xx}^{}$. In particular,
we have 
\begin{align}\label{1D-system++}
  \frac\varrho K\DDDT u-\frac{\Kv}K\DDT u_{xx}+\frac{\DT\pi_x^{}}K=0\ \ \ \text{ and }\ \ \
   \frac\varrho H\DDT u_{xx}^{}-\frac{\Kv}H\DT u_{xxxx}^{}+\frac{\pi_{xxx}^{}}H=0\,.
\end{align}
When subtracting them, we can eliminate $\pi$ by using
$\frac1K\DT\pi_x^{}-\frac1H\pi_{xxx}^{}=-\DT u_{xx}^{}$
and obtain a dispersive equation
$\varrho \DDDT u-\varepsilon\varrho\DDT u_{xx}^{}
  -\Kv \DDT u_{xx}
  -K\DT u_{xx}^{}+\varepsilon \Kv\DT u_{xxxx}^{}=0$, i.e.\ the parabolic equation
\begin{align}\label{dispersion+}
  \varrho \DDT u-(\Kv{+}\varepsilon\varrho)
 \DT u_{xx}
  -K u_{xx}^{}+\varepsilon \Kv u_{xxxx}^{}=0\,\ \ \ \ \text{ with }\ \varepsilon:=\frac KH
  \,.
\end{align}
Using again the ansatz $u={\rm e}^{{\rm i}(w t+x/\lambda)}$,
we have
$\varepsilon \Kv u_{xxxx}^{}=\varepsilon \Kv u/\lambda^4$,
and the one-dimensional 
dispersive wave equation \eq{dispersion+} yields
the algebraic condition 
\begin{align}\label{KK-relation-+}
\varrho
w^2
-{\rm i}(\Kv{+}\varepsilon\varrho)\frac w{\lambda^2}
-\frac K{\lambda^2}-\frac{\varepsilon \Kv}{\lambda^4}=0
\,.
\end{align}
Considering again $w=\omega+{\rm i}\gamma$, 
the Kramers-Kronig's relations \eq{KK-relation} now looks as
\begin{align}\label{KK-relation++}
  &\varrho(\omega^2-\gamma^2)=\frac K{\lambda^2}
  +\frac{\varepsilon \Kv}{\lambda^4}
  -(\Kv{+}\varepsilon\varrho)\frac\gamma{\lambda^2}
\ \ \ \text{ and }\ \ \
2\varrho\gamma=\frac{\Kv{+}\varepsilon\varrho}{\lambda^2}\,.
\end{align}
From the latter equation we can read that
$\gamma=(\Kv{+}\varepsilon\varrho)/(2\varrho\lambda^2)$
and then, realizing that the speed of waves $v=\omega^2\lambda^2$, the
former equation in \eq{KK-relation++} yields
\begin{align}\nonumber
v=v(\lambda)
=\sqrt{\frac K\varrho
+\frac{\varepsilon \Kv}{\varrho\lambda^2}
-(\Kv{+}\varepsilon\varrho)\frac\gamma{\varrho}+\lambda^2\gamma^2}
=\sqrt{\frac K\varrho
-\frac{(\Kv{-}\varepsilon\varrho)^2}{4\varrho^2\lambda^2}}\,.
\end{align}
Like \eq{v=v(lambda)}, we can see again the normal dispersion with the
effect that high-frequency waves cannot propagate at all.
Again, $\gamma$ determines the attenuation which naturally rises
for higher frequencies. 

\begin{remark}[{\sl Multipolar fluids}]\label{rem-multipolar}
  \upshape
  Another enhancement by dissipative gradient terms can exploit the concept of
  the 2nd-grade nonsimple fluids,
devised by E.\,Fried and M.\,Gurtin \cite{FriGur05SGFT,FriGur06TBBC}
and earlier, even more generally and nonlinearly as multipolar fluids, by
J.\,Ne\v cas at al.\ \cite{BeBlNe92PBMV,NeNoSi89GSIC,NecRuz92GSIV}. More
specifically, the 2nd-grade multipolar variant of the model \eq{NS-a} with
\eq{Oskolkov} reads as
\begin{subequations}\label{NS-modif-Bernoulli}
  \begin{align}\label{NS-a-modif-Bernoulli}
    &\varrho\DT \vv+\varrho(\vv{\cdot}\nabla)\vv
    -{\rm div}\big(\bbD e(\vv)-
    \pi\bbI-{\rm div}(\bbH\nabla e(\vv))\big)=g-\frac\varrho2({\rm div}\,\vv)\vv,
\\[-.5em]&\frac1K\DT\pi+{\rm div}\,\vv=\frac1H\Delta\pi\,
\label{NS-b-modif-Bernoulli}
\end{align}\end{subequations}
with $\bbH$ a positive-definite tensor of hyper-viscous moduli. This needs
enhancement of the boundary condition \eq{BC+}, say
\begin{align}\label{BC+multi}
  \vv_{\rm n}=0,\qquad\big(\sigma\vec{n}-\divS(\bbH\nabla e(\vv)\vec{n})\big)_{\rm t}
  +\kappa\vv_{\rm t}=0,
  \quad\text{and}\quad\bbH\nabla e(\vv)\Vdots(\vec{n}\otimes\vec{n})=0
\end{align}
where $\sigma=\bbD e(\vv)-\pi\bbI$ and where
 ``$\divS$'' is the surface divergence defined as
$\divS(\cdot)=\mathrm{tr}\big(\nablaS(\cdot)\big)$ with
$\mathrm{tr}(\cdot)$ denoting the trace and $\nablaS$ denoting the
surface gradient given by $\nablaS v=(\mathbb I- \vec n{\otimes}\vec
n)\nabla v= \nabla v-\frac{\partial v}{\partial\vec{n}}\vec{n}$. This model
gives a new term $\bbH\nabla e(\vv)\Vdots
  \nabla e(\vv)$ into the dissipation
rate in the energetics \eq{energy+}, which improves the estimate
\eq{est1-Oskolkov} by replacing $H^1(\Omega;\R^d)$ by $H^2(\Omega;\R^d)$.
The analysis of the dispersion is more involved: \eq{1D-system+} now
looks as
\begin{align}\label{1D-system+Fried-Gurtin}
  &\DT u=\vv,\ \ \ \ \ \ \ \ \varrho\DT\vv-\Kv v_{xx}+\pi_x^{}+h\vv_{xxxx}=0,\ \ \ \text{ and }\ \ \ 
  \frac{\DT\pi}K+\vv_x^{}=\frac{\pi_{xx}^{}}H\,.
\end{align}
with $h>0$ related to $\bbH=[h]$ in the 1-dimensional situation, while
\eq{1D-system++} expands by the terms $h\DDT u_{xxxx}$ and  $h\DT u_{xxxxxx}$.
The dispersive equation \eq{dispersion+} thus expands as
\begin{align*}
  \varrho \DDT u-(\Kv{+}\varepsilon\varrho)
 \DT u_{xx}
 -K u_{xx}^{}+\varepsilon \Kv u_{xxxx}^{}
+hK\DDT u_{xxxx}-h\varepsilon\DT u_{xxxxxx}
=0
  \,.
\end{align*}
Moreover, one can also think about combination of this dissipative
multipolar concept with the conservative-gradient
models as \eq{NS-Delta-pi}. This would lead to the system
\begin{subequations}\label{NS-Delta-pi-multi}
  \begin{align}
    &\!\varrho\DT \vv+\varrho(\vv{\cdot}\nabla)\vv-{\rm div}
\big(\bbD e(\vv)-(\pi{-}\ell^2\Delta\pi)\bbI-{\rm div}(\bbH\nabla e(\vv))\big)
=g-\frac\varrho2({\rm div}\,\vv)\,\vv
\label{NS-Delta-pi-a-multi}
\,,
    \\\label{NS-Delta-pi-b-multi}
    &\!\DT\pi+K\,{\rm div}\,\vv\,=0\,.
\end{align}\end{subequations}
The energetics of \eq{NS-Delta-pi-multi} when completed by the boundary
conditions \eq{BC++} can be revealed by using again the
calculus \eq{calculus-pi++}.
This gives again the energy balance \eq{energy-gradient} only with $H$
replaced by $\ell^2/K$ and with $c$ replaced by $c/\ell^2$, and here in
addition with the dissipation rate $\bbH\nabla e(\vv)\Vdots\nabla e(\vv)$.
The estimate like \eq{unique-est-} below will yield uniqueness.
The strain-rate gradient versus pressure gradient
are gradient-theoretical concepts known as Aifantis
\cite{Aifa92RGLD,AskAif11GESD,AskAif09GEFW}
versus Eringen \cite{Erin83DENE}, respectively.
In \eq{NS-Delta-pi-a-multi}, we combined both concepts.
\end{remark}

\section{Truly convective ``semi-compressible'' models}\label{sec-convective}
All the above models may still be referred under the name ``quasi'' because
they are not fully mechanically consistent when mixing convective
time derivative in \eq{NS-a} or in \eq{NS-Delta-pi-a} or
\eq{NS-a-modif-Bernoulli}, i.e.\ the Eulerian description, with the partial
time derivative in \eq{NS-b}, \eq{Sobolev}, \eq{Oskolkov}, or
\eq{NS-b-modif-Bernoulli}, i.e.\ the Lagrangian description.
Now we want to improve (some of) these models to be formulated consistently
fully convectively, i.e.\ in Eulerian coordinates.
It is also likely that the criticism \cite{OtStLi09IDCM} of the diffusive term
in \eq{Oskolkov} or \eq{NS-modif-Bernoulli} because of lack of Galilean
invariance is thus suppressed so that the justification of these
diffusive terms in \cite{VaPaGr17EMFF} is fully in effect.

In fact, the derivation of the quasi-incompressible model \eq{NS-b} from
the fully-compressible equation \eq{compress-NS1} under the 
small perturbation ansatz \eq{ansatz} was not precise. 
In fact, 
from \eq{compress-NS1} when using
\eq{ansatz}, we have
\begin{align}\label{convective-}
  \varrho\frac1K\DT\pi
  +\vv{\cdot}\nabla\Big(\varrho\frac\pi K\Big)
  +\varrho\Big(1{+}\frac\pi K\Big){\rm div}\,\vv=0\,,
\ \ \ \text{ i.e. }\ \ \ 
\frac{\DT\pi{+}\vv{\cdot}\nabla\pi}{K{+}\pi}+{\rm div}\,\vv=0\,.
\end{align}
Instead of \eq{NS-b}, it should lead rather to
\begin{align}\label{convective}
 \frac1K\big(\DT\pi{+}\vv{\cdot}\nabla\pi\big)+{\rm div}\,\vv=0
\end{align}
when assuming $|\pi|\ll K$ to replace the coefficient $1/(K{+}\pi)$ by $1/K$.
Mechanically, one can say that replacing the extensive variable (density)
$\rho$  in \eq{compress-NS} by the intensive variable (pressure) $\pi$
in \eq{NS}, should lead to replacement the transport equation \eq{compress-NS1}
for the extensive variable $\rho$  by a transport equation
\eq{convective} for the intensive variable $\pi$;
  let us remind that intensive variables are those whose magnitude is
  independent of the size of the system (volume) whereas the extensive ones
  are those whose magnitude is additive for subsystems.

  However, replacement of \eq{compress-NS1} by \eq{convective}
  brings analytical difficulties 
  because $\nabla\pi$, similarly as in the nonconvective variant  \eq{NS-a},
  is not a-priori estimated in the model \eq{NS-a} with  \eq{convective}.
  Even the by-part integration $\int_\Omega(\vv{\cdot}\nabla\pi)\wt\pi\,\d x=
  -\int_\Omega\pi(\vv{\cdot}\nabla\wt\pi)+\pi({\rm div}\,\vv)\wt\pi\,\d x$
  which could eliminate $\nabla\pi$ from a
  weak formulation seems problematic because then one would need a strong
  convergence in $\nabla\vv$ or in $\pi$, which does not seem directly
  at disposal. Similarly, the convective modification of the model \eq{NS-a}
  with \eq{Sobolev} from Section \ref{sec-dispersive} seems analytically
  problematic.

Yet, some of the above discussed models bear the physically relevant
convective modification which is simultaneously amenable for analysis.
In particular, for \eq{Oskolkov}, the force-equilibrium
equation \eq{NS-a} is to be now slightly modified
to obtain the model
\begin{subequations}\label{Oskolkov-convective}
  \begin{align}
    &\varrho\DT\vv+\varrho(\vv{\cdot}\nabla)\vv
    -{\rm div}\Big(\bbD e(\vv)
    -\Big(\frac{\pi^2}{2K}+\pi\Big)\bbI\Big)
  =g-\frac\varrho2({\rm div}\,\vv)\,\vv\,,
 \label{NS-modif++}
\\[-.2em]
&\label{Oskolkov+}
\frac1K(\DT\pi{+}\vv{\cdot}\nabla\pi)+{\rm div}\,\vv=\frac1H\Delta\pi\,.
\end{align}\end{subequations}
Beside the 
``hydrostatic''  pressure $\pi$, we now also see an  ``internal'' 
pressure due to the elastic energy $\pi^2/(2K)$ in \eq{NS-modif++}.
 This contribution due to internal energy (here elastic but it can also
be e.g.\ magnetic or chemical) is a thermodynamically justified effect,
cf.\ e.g.\ \cite{Mark13IPLS}, which disappears only in ideally
incompressible models.
Here, it  arises when using the Green formula
\begin{align}\nonumber
  \int_\Omega\pi\vv{\cdot}\nabla\pi\,\d x
&=\int_\Gamma\pi^2\vv\cdot\vec{n}\,\d S-\int_\Omega{\rm div}(\pi\vv)\pi\,\d x
  \\\nonumber
  &=\int_\Gamma\pi^2\vv\cdot\vec{n}\,\d S
  -\int_\Omega(\nabla\pi\cdot\vv)\pi+\pi^2{\rm div}\,\vv\,\d x
  \\&=\int_\Gamma\frac12\pi^2\vv\cdot\vec{n}\,\d S
  -\int_\Omega\frac12\pi^2{\rm div}\,\vv\,\d x,
\end{align}
so that the calculus \eq{calculus-pi+++} now enhances as
\begin{align}\nonumber
  &\int_\Omega\nabla\pi\cdot \vv\,\d x-\int_\Gamma\pi\,(\vv\cdot\vec{n})\,\d S
=\int_\Omega\pi\Big(\frac1K\DT\pi+\frac1K\vv{\cdot}\nabla\pi
  -\frac1H\Delta\pi\Big)\,\d x
    \\&=\frac{\d}{\d t}\int_\Omega\frac{\pi^2}{2K}\,\d x
    +\int_\Omega\frac1H|\nabla\pi|^2\!-\frac{\pi^2}{2K}{\rm div}\,\vv\,\d x
    +\int_\Gamma\frac{\pi^2}{2K}\vv{\cdot}\vec{n}
    -\frac1H\pi(\nabla\pi{\cdot}\vec{n})\,\d S\,.
    \label{calculus-pi++++}\end{align}
The boundary conditions can again be considered as in \eq{BC+}.

The energetics \eq{energy+} remains the same because the additional
terms $\frac1{2K}\pi^2$ in \eq{calculus-pi++++} just cancel with the
additional pressure terms coming by the test of \eq{NS-modif++}
by $\vv$ and using the enhanced boundary condition \eq{BC++}.

\begin{proposition}[Existence of weak solutions to \eq{Oskolkov-convective}]\label{prop-3}
  Let $d\le3$, $\varrho,K,H,c>0$, $\kappa\ge0$, $\bbD$ be symmetric
  positive definite, and \eq{ass} hold. Then the system
  \eq{Oskolkov-convective} with the boundary conditions \eq{BC+}
  possesses a weak solution 
$\vv\in C_{\rm w}(I;L^2(\Omega;\R^d))\:\cap\:L^2(I;H^1(\Omega;\R^d))$
  and $\pi\in C_{\rm w}(I;H^1(\Omega))$.
\end{proposition}

\noindent{\it  Sketch of the proof.}
  Again we consider the Galerkin approximation.
The a-priori estimates 
 (\ref{est-Oskolkov}a,b)  remain
the same. The a-priori estimate for $\varrho\DT\vv_k$ 
(in its Hahn-Banach extension) in 
in $L^2(I;H^1(\Omega;\R^d)^*)+L^{5/4}(I{\times}\Omega;\R^d)$
and the estimate \eq{est6-Oskolkov} are slightly
weakened here because of the mentioned additional
  pressure $\frac1{2K}\pi^2$ in \eq{NS-modif++}, bounded in
$L^1(I;L^3(\Omega))$, and because of the
term $\vv{\cdot}\nabla\pi$ in \eq{Oskolkov+}.
As for $\DT\pi$, in fact we can now only estimate it by comparison
$\DT\pi_k=\frac KH\Delta\pi_k
-K{\rm div}\,\vv_k-\vv_k{\cdot}\nabla\pi_k$ understood in its Galerkin
identity after its Hahn-Banach extension uniformly in the
space $L^\infty(I;H^1(\Omega)^*)+L^2(I;L^{3/2}(\Omega))$
when realizing that $\vv{\cdot}\nabla\pi\in L^2(I;L^{3/2}(\Omega))$
for $d=3$. 
Thus we can rely on
\begin{subequations}\label{est+-Oskolkov-modif}\begin{align}
    &\label{est5-Oskolkov-modif}
    \|\varrho\DT \vv_k\|
    _{L^1(I;H^1(\Omega;\R^d)^*)+L^{5/4}(I{\times}\Omega;\R^d)}^{}\le C\ \ \text{ and }\ \ 
    \\&\label{DTpi-est}
    \|\DT\pi_k\|_{L^\infty(I;H^1(\Omega)^*)+L^2(I;L^{3/2}(\Omega))}^{}\le C\,.
\end{align}\end{subequations}

The convergence
in the additional nonlinear term
$\vv_k\cdot\nabla\pi_k$ in \eq{Oskolkov+} is easy because $\nabla\pi_k$
converge weakly* in $L^2(I{\times}\Omega;\R^d)$ due to \eq{est2-Oskolkov} while
$\vv_k$ converge strongly in $L^2(I{\times}\Omega;\R^d)$ due to
\eq{est1-Oskolkov} and the Aubin-Lions compact-embedding theorem using
also \eq{est5-Oskolkov-modif}.
In addition, we need strong convergence of $\pi_k$ in the additional
pressure $\frac1{2K}\pi_k^2$, which is again simple when using the Aubin-Lions
theorem with \eq{DTpi-est} with the only technicality that
$\DT\pi_k$ in \eq{DTpi-est} is to be understood as a Hahn-Banach extension
or the estimate \eq{DTpi-est} is to be weakened by using only seminorms,
cf.\ \cite[Sect.8.4]{Roub13NPDE}.
\hfill$\Box$

\medskip

  
One should note that the rigorous energy conservation behind the model
\eq{Oskolkov-convective} is not granted by our estimates because
the inertial force $\varrho\DT\vv$ does not enjoy enough regularity
to be eligible for being tested by $\vv$. Also the uniqueness
for $d=3$ is not obvious, cf.\ Remark~\ref{rem-uni}. For these reasons, one
may still consider
the fully convective variant of the multipolar model \eq{NS-modif-Bernoulli},
which looks as
  \begin{subequations}\label{NS-modif-Bernoulli-convective}
  \begin{align}\label{NS-a-modif-Bernoulli-convective}
    &\varrho\DT \vv+\varrho(\vv{\cdot}\nabla)\vv
    -{\rm div}\Big(\bbD e(\vv)
    -\Big(
    \frac{\pi^2}{2K}{+}\pi\Big)\bbI
    -{\rm div}(\bbH\nabla e(\vv))\Big)=g-\frac\varrho2({\rm div}\,\vv)\vv,
  \\[-.3em]&\label{NS-b-modif-Bernoulli-convective}
  \frac1K\big(\DT\pi{+}\vv{\cdot}\nabla\pi\big)+{\rm div}\,\vv
  =\frac1H\Delta\pi\,.
  \end{align}\end{subequations}
  {}together with the boundary conditions  \eq{BC+}
  combined with \eq{BC+multi}, i.e.\ here
  \begin{align}\nonumber
 & \vv_{\rm n}=0,\qquad\big(\sigma\vec{n}-\divS(\bbH\nabla e(\vv)\vec{n})\big)_{\rm t}
    +\kappa\vv_{\rm t}=0\ \ \text{ with }\ \
\sigma=\bbD e(\vv)-\Big(\frac{\pi^2}{2K}{+}\pi\Big)\bbI\,,
\\[-.3em]\label{BC+multi+}
  &\bbH\nabla e(\vv)\Vdots(\vec{n}\otimes\vec{n})=0,\ \ \text{ and }\ \
 \frac1{H}\nabla\pi\cdot\vec{n}+\frac\pi c=0\,.
  \end{align}  
This allows for improving the estimate of $\DT\pi$ from \eq{DTpi-est}.
Another important attribute of the semi-compressible convective multipolar model
\eq{NS-modif-Bernoulli-convective} is that it allows for uniqueness even in
3-dimensional situations.

\begin{proposition}[Weak solutions to
    \eq{NS-modif-Bernoulli-convective}--\eq{BC+multi+}]\label{prop-4}
  Let  $d\le3$, $\varrho,K,H,c>0$, $\kappa\ge0$, 
  $\bbD$ and $\bbH$ be symmetric
  positive definite, and let 
  \begin{align}\label{ass+}
 \!\!g\!\in\! L^2(I{\times}\Omega;\R^d),\ \ \ 
  \vv_0\!\in\!H^1(\Omega;\R^d),\ \ \text{ and }\ \pi_0\!\in\!H^1(\Omega)
\end{align}
hold. Then:\\
\Item{(i)}{the system \eq{NS-modif-Bernoulli-convective} with the initial and
  boundary conditions \eq{IC} and \eq{BC+multi}
  has a unique  weak solution
  $\vv\in C_{\rm w}(I;H^1(\Omega;\R^d))\:\cap\:L^2(I;H^2(\Omega;\R^d))$
  and $\pi\in C_{\rm w}(I;H^1(\Omega))\:\cap\:H^1(I;L^2(\Omega))$.
Moreover, $\Delta\pi\in L^2(I{\times}\Omega)$.}
\Item{(ii)}{This weak solution conserves energy. More specifically,
    for any $t\in I$ it holds}
    \begin{align}\nonumber
&\int_\Omega\frac\varrho2|\vv(t)|^2+
 \frac{\pi(t)^2}{2K}\,\d x
 +\int_0^t\!\!\int_\Omega\bbD e(\vv){:}e(\vv)
 +\bbH\nabla e(\vv)\Vdots\nabla e(\vv)+\frac{|\nabla\pi|^2}H\,\d x\d t
  \\[-.1em]&\qquad\qquad+\int_0^t\!\!\int_\Gamma
    \kappa|\vv_{\rm t}|^2+\frac{\pi^2}c\,\d S\d t 
  =\int_\Omega\frac\varrho2|\vv_0|^2+\frac{\pi_0^2}{2K}
  \,\d x+\int_0^t\!\!\int_\Omega g{\cdot}\vv\,\d x\d t
  \,.
\label{energy++}\end{align}  
  \Item{(iii)}{If $\Omega$ is smooth, then even
    $\vv\in L^2(I;H^4(\Omega;\R^d))$ and 
  $\pi\in L^2(I;H^2(\Omega))$
    and \eq{NS-modif-Bernoulli-convective} holds a.e.\ on $I{\times}\Omega$.
}
\end{proposition}

  \noindent{\it Proof.}
  First, testing \eq{NS-modif++} in its Galerkin approximation by the Galerkin
  approximation $\vv_k$ of $\vv$ and using the calculus \eq{calculus-pi++++}
  for \eq{NS-b-modif-Bernoulli-convective}, we obtain the energy balance \eq{energy++} for the
  approximate solutions. Of course, \eq{NS-b-modif-Bernoulli-convective} is
  understood here in its Galerkin approximation for which the test by
$\pi_k$ used in \eq{calculus-pi++++} is legitimate.
  From this, one can read the a-priori estimates 
  \begin{subequations}\label{est-Oskolkov-nonsimple}
    \begin{align}\label{est1-Oskolkov-nonsimple}
    &\|\vv_k\|_{L^\infty(I;L^2(\Omega;\R^d))\:\cap\:L^2(I;H^2(\Omega;\R^d))}^{}\le C\,,
    \\&\label{est2-Oskolkov-nonsimple}\|\pi_k\|_{L^\infty(I;L^2(\Omega))\:\cap\:L^2(I;H^1(\Omega))}^{}\le C\,.
\end{align}\end{subequations}

  Second, by testing \eq{NS-b-modif-Bernoulli-convective}
by $\DT\pi_k$, we obtain
\begin{align*}
 &\int_\Omega\frac1K\big|\DT\pi_k\big|^2
  +\frac1{2H}\frac{\pl}{\pl t}|\nabla\pi_k|^2\,\d x
+\int_\Gamma\frac1{2c}\frac{\pl}{\pl t}|\pi_k|^2\,\d S
=
  -\int_\Omega\Big(\frac1K\vv_k{\cdot}\nabla\pi_k+{\rm div}\,\vv_k\Big)\DT\pi_k\,\d x
  \\[-.2em]&\qquad\qquad
  \le\frac1{K}\|\vv_k\|_{L^\infty(\Omega;\R^d)}^2\|\nabla\pi\|_{L^2(\Omega;\R^d)}^2
  +K\|{\rm div}\,\vv_k\|_{L^2(\Omega)}^2
  +\frac1{2K}\|\DT\pi_k\|_{L^2(\Omega)}^2\,.
\end{align*}
Using the embedding $H^2(\Omega)\subset L^\infty(\Omega)$ so that
\begin{align}\label{for-Gronwall}
  \Big\{t\mapsto\|\vv_k(t)\|_{L^\infty(\Omega;\R^d)}^2\Big\}_{k\in\N}\
  \text{ is bounded in }\  L^1(I)
\end{align}
due to
\eq{est1-Oskolkov-nonsimple}, by the Gronwall inequality we improve the
estimate \eq{est2-Oskolkov-nonsimple} as
\begin{align}\label{est3-Oskolkov-nonsimple}
  \|\pi_k\|_{L^\infty(I;H^1(\Omega))\:\cap\:H^1(I;L^2(\Omega))}^{}\le C\,.
\end{align}

Third, let us realize that,
by \eq{est3-Oskolkov-nonsimple},
${\rm div}(\pi_k^2\bbI/2)=\pi_k\nabla\pi_k$ is bounded in 
$L^\infty(I;L^3(\Omega;\R^d))$. Therefore,
we can test \eq{NS-a-modif-Bernoulli-convective}
in its Galerkin approximation by $\DT v_k$. Thus
\begin{align*}
  &\int_\Omega\varrho|\DT\vv_k|^2\,\d x
  +\frac12\frac{\d}{\d t}\int_\Omega\bbD e(\vv_k):e(\vv_k)
  +\bbH\nabla e(\vv_k):\nabla e(\vv_k)\,\d x=
  \\&=\int_\Omega g\cdot\DT\vv_k
  -\frac{\pi_k}{2K}\nabla\pi_k\cdot\DT\vv_k-\nabla\pi_k\cdot\DT\vv_k
  -(\vv_k\cdot\nabla)\vv_k\cdot\DT\vv_k
  -\frac\varrho2({\rm div}\,\vv_k)\vv_k\cdot\DT\vv_k\,\d x\,.
\end{align*}
{}The first right-hand side term is to be estimated as
$\int_\Omega g\cdot\DT\vv_k\,\d x\le\frac\varrho2\|\DT\vv_k\|_{L^2(\Omega;\R^d)}^2+
+\frac1{2\varrho}\|g\|_{L^2(\Omega;\R^d)}^2$; here the assumption
$g\in L^2(I{\times}\Omega;\R^d)$ is employed. 
When estimating the last two term as
\begin{align*}
  -\int_\Omega(\vv_k{\cdot}\nabla)\vv_k\cdot\DT\vv_k+
  \frac\varrho2({\rm div}\,\vv_k)\vv_k\cdot\DT\vv_k\,\d x
  \le \frac\varrho2\|\DT\vv_k\|_{L^2(\Omega;\R^d)}^2\qquad\qquad\qquad
  \\
  +\ C_\varrho^{}\|\vv_k\|_{L^\infty(\Omega;\R^d)}^2\big(\|\vv_k\|_{L^2(\Omega;\R^d)}^2
+\|e(\vv_k)\|_{L^2(\Omega;\R^{d\times d})}^2\big)
\end{align*}
with some $C_\varrho^{}$ depending on $\varrho$ and involving also a
constant from the Korn inequality on $\Omega$, and when using
again \eq{for-Gronwall}, by the Gronwall inequality 
we still obtain the estimate
\begin{align}\label{est4-Oskolkov-nonsimple}
\|\DT\vv_k\|_{L^2(I{\times}\Omega;\R^d)}^{}\le C\ \ \ \text{ and }\ \ \ 
\|e(\vv_k)\|_{L^\infty(I;H^1(\Omega;\R^{d\times d}))}^{}\le C\,.
\end{align}

The limit passage of selected weakly* convergent subsequences towards weak
solutions $(\vv,\pi)$ to \eq{NS-modif-Bernoulli-convective} is then easy. 
Then, from \eq{NS-b-modif-Bernoulli-convective} we can also see that
$\Delta\pi=\frac HK(\DT\pi+\vv{\cdot}\nabla\pi)+H{\rm div}\,\vv\in
L^2(I{\times}\Omega)$ due to \eq{est3-Oskolkov-nonsimple}.

\COMMENT{We have thus ${\rm div}\,\vv\in L^\infty(I;L^6(\Omega))$ and
  $\vv{\cdot}\nabla\pi\in L^\infty(I;L^2(\Omega))$.
  DO NOT WE HAVE then $\pi\in L^\infty(I;H^2(\Omega))$ from the
  linear parabolic equation $\frac1K\DT\pi-\frac1H\Delta\pi\in
  =-{\rm div}\,\vv-\frac1K\vv{\cdot}\nabla\pi\in L^\infty(I;L^2(\Omega))$. ?????}

As to the uniqueness,
enhancing the calculus \eq{NS-uniqueness+}, we have for a.a.\ time
instances $t\In I$ (with $t$ omitted in the following formulas for notational
simplicity) that
\begin{align}
\nonumber&\!\!\!
\|\sqrt\bbD e(\vv_{12})\|_{L^2(\Omega;\R^{d\times d})}^2
+\|\sqrt\bbH\nabla e(\vv_{12})\|_{L^2(\Omega;\R^{d\times d\times d})}^2
\\\nonumber&\hspace{5em}+\frac1{2H}\|\nabla\pi_{12}\|_{L^2(\Omega;\R^d)}^2
+\frac{\d}{\d t}\Big(\frac\varrho2\|\vv_{12}\|_{L^2(\Omega;\R^d)}^2
+\frac1{2K}\|\pi_{12}\|_{L^2(\Omega)}^2\Big)
\\\nonumber&\hspace{1em}
=\int_\Omega\varrho\big((\vv_2{\cdot}\nabla)\vv_2-\vv_1{\cdot}\nabla)\vv_1\big)
\cdot\vv_{12}
+\frac\varrho2\big(({\rm div}\,\vv_2)\vv_2-({\rm div}\,\vv_1)\vv_1\big)
\cdot\vv_{12}
\\[-.5em]\nonumber&\hspace{5em}
+\frac1K\big(|\pi_1|^2-|\pi_2|^2\big)\bbI:e(\vv_{12})
+\frac1K\big(\vv_1{\cdot}\nabla\pi_1-\vv_2{\cdot}\nabla\pi_2\big)\pi_{12}
\,\d x\!
\\\nonumber&\hspace{1em}
=\int_\Omega\varrho\Big((\vv_{12}{\cdot}\nabla)\vv_1+\frac12({\rm div}\,\vv_{12})\,\vv_2
+(\vv_2{\cdot}\nabla)\vv_{12}+\frac12({\rm div}\,\vv_1)\,\vv_{12}\Big)\cdot \vv_{12}
\\[-.5em]&\hspace{5em}
+\frac1K\Big((\pi_1{+}\pi_2)\,{\rm div}\,\vv_{12}
+\vv_{12}{\cdot}\nabla\pi_1+\vv_2{\cdot}\nabla\pi_{12}
\Big)\pi_{12}\,\d x\,.
\label{unique-est-}\end{align}
where $\vv_{12}:=\vv_1{-}\vv_2$ and $\pi_{12}:=\pi_1{-}\pi_2$ for two weak
solutions $(\vv_1,\pi_1)$ and $(\vv_2,\pi_2)$.
The particular terms on the right-hand side can be estimated by H\"older's and
Young's inequalities. In particular, 
by the Gagliadro-Nirenberg inequality together with Korn's inequality
 for the 2nd-grade nonsimple materials \cite[Sect.\,5.2]{KruRou19MMCM},
we have 
\begin{align*}
\|\vv_{12}\|_{L^r(\Omega;\R^d)}^2&\le C_r^{}\Big(
\|\vv_{12}\|_{L^2(\Omega;\R^d)}^{}\|\nabla^2\vv_{12}\|_{L^2(\Omega;\R^{d\times d\times d})}^{}
+\|\vv_{12}\|_{L^2(\Omega;\R^d)}^2\Big)
\\
&\le C_r'\Big(
\|\vv_{12}\|_{L^2(\Omega;\R^d)}^{}\|\nabla e(\vv_{12})\|_{L^2(\Omega;\R^{d\times d\times d})}^{}
+\|\vv_{12}\|_{L^2(\Omega;\R^d)}^2\Big)
\end{align*} 
for any $1\le r\le 6$.
Here we will use it for $r=12/5$. Thus, we can estimate
\begin{subequations}\label{unique-est}\begin{align}
  &\nonumber\int_\Omega\varrho\Big((\vv_{12}{\cdot}\nabla)\vv_1
  +\frac12({\rm div}\,\vv_1)\,\vv_{12}\Big)\cdot \vv_{12}\,\d x
  \le\frac{3\varrho}2\|e(\vv_1)\|_{L^6(\Omega;\R^{d\times d})}\|\vv_{12}\|_{L^{12/5}(\Omega;\R^d)}^2
  \\&\quad
  \nonumber\le
  \frac{3\varrho C_r'}2\|e(\vv_1)\|_{L^6(\Omega;\R^{d\times d})}
  \|\vv_{12}\|_{L^2(\Omega;\R^d)}
  \Big(\|\vv_{12}\|_{L^2(\Omega;\R^d)}\!+
  \|\nabla e(\vv_{12})\|_{L^2(\Omega;\R^{d\times d\times d})}\Big)
  \\&\nonumber\qquad
  \le\epsilon\|\nabla e(\vv_{12})\|_{L^2(\Omega;\R^{d\times d\times d})}^2
\\&\qquad\
+\Big(\frac{3\varrho C_r'}2\|e(\vv_1)\|_{L^6(\Omega;\R^{d\times d})}^{}
+
\frac{9\varrho^2{C_r'}^2}{16\epsilon}
  \|e(\vv_1)\|_{L^6(\Omega;\R^{d\times d})}^2
\Big)
  \|\vv_{12}\|_{L^2(\Omega;\R^d)}^2,\!
\\  \nonumber&\int_\Omega\varrho\Big(\frac12({\rm div}\,\vv_{12})\,\vv_2
  +(\vv_2{\cdot}\nabla)\vv_{12}\Big)\cdot \vv_{12}\,\d x
  \\&\qquad\qquad\le\nonumber
  \frac{3\varrho}2\|\vv_2\|_{L^\infty(\Omega;\R^d)}
  \|e(\vv_{12})\|_{L^2(\Omega;\R^{d\times d})}\|\vv_{12}\|_{L^2(\Omega;\R^d)}
  \\&\qquad\qquad\le\epsilon\|e(\vv_{12})\|_{L^2(\Omega;\R^{d\times d})}^2
  +\frac{9\varrho^2}{16\epsilon} \|\vv_2\|_{L^\infty(\Omega;\R^d)}^2
 \|\vv_{12}\|_{L^2(\Omega;\R^d)}^2\,,
 \\\nonumber&\int_\Omega\frac1K(\pi_1{+}\pi_2)({\rm div}\,\vv_{12})\pi_{12}\,\d x
 \le\frac1K\|\pi_1{+}\pi_2\|_{L^4(\Omega)}\|e(\vv_{12})\|_{L^4(\Omega;\R^{d\times d})}
 \|\pi_{12}\|_{L^2(\Omega)}
 \\&\qquad\qquad\le\epsilon \|e(\vv_{12})\|_{L^4(\Omega;\R^{d\times d})}^2
  +\frac1{4\epsilon K^2}\|\pi_1{+}\pi_2\|_{L^4(\Omega)}^2\|\pi_{12}\|_{L^2(\Omega)}^2\,,
\\  \nonumber&\int_\Omega\frac1K
(\vv_{12}{\cdot}\nabla\pi_1)\pi_{12}\,\d x\le\frac1K
\|\vv_{12}\|_{L^\infty(\Omega;\R^d)}\|\nabla\pi_1\|_{L^2(\Omega;\R^d)}\|\pi_{12}\|_{L^2(\Omega)}
\\&\qquad\qquad\le\epsilon\|\vv_{12}\|_{L^\infty(\Omega;\R^d)}^2
+\frac1{4\epsilon K^2}\|\nabla\pi_1\|_{L^2(\Omega;\R^d)}^2\|\pi_{12}\|_{L^2(\Omega)}^2\,,
\ \ \text{ and}
\\\nonumber&\int_\Omega\frac1K
  (\vv_2{\cdot}\nabla\pi_{12})\pi_{12}\,\d x\le\|\vv_2\|_{L^\infty(\Omega;\R^d)}
  \|\nabla\pi_{12}\|_{L^2(\Omega;\R^d)}\|\pi_{12}\|_{L^2(\Omega)}
  \\&\qquad\qquad\le\epsilon \|\nabla\pi_{12}\|_{L^2(\Omega;\R^d)}^2
  +\frac1{4\epsilon K^2}\|\vv_2\|_{L^\infty(\Omega;\R^d)}^2\|\pi_{12}\|_{L^2(\Omega)}^2\,.
\end{align}\end{subequations}
Taking $\epsilon>0$ small enough, the $\epsilon$-terms on the
right-hand sides of \eq{unique-est} can be absorbed in the left-hand
side of \eq{unique-est-} and then, using that
$t\mapsto\|e(\vv_{1}(t))\|_{L^6(\Omega;\R^{d\times d})}^2$, 
$t\mapsto\|\pi_1(t){+}\pi_2(t)\|_{L^4(\Omega)}^2$ 
and 
$t\mapsto\|\vv_{1}(t)\|_{L^\infty(\Omega;\R^d)}^2+\|\nabla\pi_1(t)\|_{L^2(\Omega;\R^d)}^2$
are
$L^1(I)$, we treat it by Gronwall inequality.

As to (ii), let us realize that, by \eq{est4-Oskolkov-nonsimple},
$\varrho\DT\vv\in L^2(I;H^2(\Omega;\R^d)^*)$ is certainly in duality with
$\vv\in L^2(I;H^2(\Omega;\R^d))$.
Moreover, $\DT\pi$, $\vv{\cdot}\nabla\pi$, and $\Delta\pi$ belong to
$L^2(I{\times}\Omega)$ so that these terms allows for testing by $\pi$, as
used in the calculus \eq{calculus-pi++++}.
Also ${\rm div}^2(\bbH\nabla e(\vv))\in L^2(I;H^2(\Omega;\R^d)^*)$
is in duality with (and allows for the test by)
$\vv\in L^2(I;H^2(\Omega;\R^d))$. Thus the
calculus \eq{calculus-pi++++} holds rigorously and the test
of \eq{NS-a-modif-Bernoulli-convective} leading to the energy equality
\eq{energy++} is legitimate.

Eventually, (iii) follows by the $H^2$-regularity of the $\Delta$-operator
on smooth domains the left-hand side of \eq{NS-b-modif-Bernoulli-convective}
is valued in $L^2(\Omega)$. Note that $\DT\pi\in L^2(I;L^2(\Omega))$ due to
\eq{est3-Oskolkov-nonsimple} while $\vv{\cdot}\nabla\pi$ and ${\rm div}\,\vv$
are in a smaller space due to \eq{est-Oskolkov-nonsimple}.
Similarly, we can use the $H^4$-regularity for the bi-harmonic operator
${\rm div}^2(\bbH\nabla e(\vv))$ and that $\DT v\in L^2(I;L^2(\Omega;\R^d))$
due to \eq{est4-Oskolkov-nonsimple}, and to obtain the
$L^2(I;H^4(\Omega;\R^d))$ regularity of $\vv$ from
\eq{NS-a-modif-Bernoulli-convective}.
\hfill$\Box$

\medskip

\begin{remark}[{\sl Pressure-dependent compressibility II}]\label{rem-K=K(pi)+}
  \upshape
  The convective modification of the incompressibility \eq{convective}
is actually sometimes considered in magnetohydrodynamics of liquid metals, cf.\
e.g.\ \cite[Formula (2.14)]{Bisk93NMHD} or \cite[Formula (9.4)]{Schn09LMHD},
where $K$ is considered as a linear function of $\pi$.
\end{remark}

\begin{remark}[{\sl Bernoulli principle}]
  \upshape
  Using the calculus $(\vv{\cdot}\nabla)\vv=(\nabla\times\vv)\times\vv+\frac12\nabla|\vv|^2$, the model \eq{Oskolkov-convective} can be written in the form
\begin{subequations}\label{Oskolkov-convective+}
  \begin{align}
    &\varrho\DT\vv+\!\!\!\!\!\!\linesunder{\varrho
(\nabla{\times}\vv){\times}\vv}{Lamb}{force}\!\!\!\!\!\!
    -{\rm div}\Big(\bbD e(\vv)
    -\Big(\!\!\!\!\!\linesunder{\frac\varrho2|\vv|^2\!+\frac{\pi^2}{2K}+\pi}{``Bernoulli''}{pressure}\!\!\!\!\!\Big)\bbI\Big)
  =g-\!\!\!\!\!\linesunder{\frac\varrho2({\rm div}\,\vv)\,\vv}{Temam's}{force}\!\!\!\!\!,
 \label{NS-modif+++}
\\[-.7em]
&\label{Oskolkov+Bernoulli}
\frac1K(\DT\pi{+}\vv{\cdot}\nabla\pi)+{\rm div}\,\vv=\frac1H\Delta\pi\,.
\end{align}\end{subequations}
This recovers the celebrated {\it Daniel Bernoulli's principle}
  \cite{Bern38HSVM} which
states that an increase in the speed of a fluid occurs
simultaneously with a decrease in pressure and vice versa.
The ``Bernoulli'' pressure in \eq{Oskolkov-convective+}
is thus the sum of the so-called dynamic pressure $\frac\varrho2|\vv|^2$,
the pressure related with the internal energy of the fluid $\frac1{2K}\pi^2$, 
and the ``hydrostatic'' pressure $\pi$.
The model \eq{NS-modif-Bernoulli-convective}
bears an analogous form and is thus
consistent with the Bernoulli principle, too.
  \end{remark}

  \begin{remark}[{\sl Convective variant of \eq{NS-Delta-pi-multi}}]
    \upshape
    The model \eq{NS-Delta-pi-multi} in its fully convective variant
  would interestingly involve also the Korteweg-like stress
$\frac{\ell^2}K(\nabla\pi\otimes\nabla\pi-\frac12|\nabla\pi|^2\bbI)$,
which arises from the calculus \eq{Korteweg} below.
The resulting system is:
\begin{subequations}\label{NS-Delta-pi-convect}
  \begin{align}\nonumber
    &\varrho\DT \vv+\varrho(\vv{\cdot}\nabla)\vv-{\rm div}
    \Big(\bbD e(\vv)-
    \Big(
    \frac{\pi^2}{2K}{+}\ell^2\frac{|\nabla\pi|^2}{2K}
         {+}\pi{-}\ell^2\Delta\pi\Big)\bbI
    \\&\hspace{9em}+\frac{\ell^2}K\nabla\pi\otimes\nabla\pi
    -{\rm div}(\bbH\nabla e(\vv))\Big)
=g-\frac\varrho2({\rm div}\,\vv)\,\vv
\label{NS-Delta-pi-a-convect}
\,,
    \\\label{NS-Delta-pi-b-convect}
    &\DT\pi+\vv{\cdot}\nabla\pi+K\,{\rm div}\,\vv\,=0\,.
\end{align}\end{subequations}
The previous calculus \eq{calculus-pi++} now enhances as
(for $\vv\cdot\vec{n}=0$ and $\nabla\pi\cdot\vec{n}=0$\COMMENT{???})
\begin{align}\nonumber
  &\int_\Omega\nabla(\pi{-}\ell^2\Delta\pi)\cdot \vv\,\d x
  =
  \int_\Omega(\pi{-}\ell^2\Delta\pi)\Big(\frac1K\DT\pi+\frac1K\vv{\cdot}
    \nabla\pi
  \Big)\,\d x
  \\[-.2em]&\nonumber
  \ \ =\frac{\d}{\d t}\int_\Omega\frac{\pi^2}{2K}
  +\frac{\ell^2|\nabla\pi|^2}{2}\,\d x
  +\int_\Omega\frac{\ell^2} K(\vv{\cdot}\nabla\pi)\Delta\pi
  -\frac{\pi^2}{2K}{\rm div}\,\vv\,\d x
    \end{align}
and then to use the calculus (by Green's formula)
\begin{align}\label{Korteweg}
  \int_\Omega(\vv{\cdot}\nabla\pi)\Delta\pi\,\d x
  &=\int_\Omega\frac12|\nabla\pi|^2{\rm div}\,\vv
  -(\nabla\pi\otimes\nabla\pi):e(\vv)\,\d x\,,
\end{align}
cf.\ e.g.\ \cite[Formula (7.7.27)]{KruRou19MMCM}.
This gives the same energetics as the model \eq{NS-Delta-pi-multi}.
However, the limit passage in the Galerkin approximation would now need also
the strong convergence of $\nabla\pi$ due to the Korteweg stress in
\eq{NS-Delta-pi-a-convect}, which seems difficult.
In principle, a higher-gradient ``multipolar-type'' term like
$-{\rm div}^2\nabla^2\pi$ instead of 0 in
\eq{NS-Delta-pi-b-convect} with correspondingly enhanced boundary conditions
would help.  A certain variant of this model even with $\varrho$ as
a variable subjected to the continuity equation, combined
also with the diffusion in \eq{NS-b-modif-Bernoulli-convective} and
amenable for existence of solutions, has been developed in \cite{BuFeMa19CCVR}.
For a different, incompressible variant with $\varrho$ fixed 
combining incompressible model with a convective-diffusive transport
of a scalar pressure-like variable, amenable for existence of solutions,
see \cite{BMPS18PDEA}.  
\end{remark}

\begin{remark}[{\sl Pressure-dependent compressibility III}]\label{rem-K=K(pi)++}
  \upshape
  As in Remark~\ref{rem-K=K(pi)}, $K=K(\pi)$ can be considered, too.
Now, \eq{convective-} suggests to take it so that $K'(0)=1$, cf.\ also
Remark~\ref{rem-K=K(pi)+}. We thus
imitate the factor $1/(K{+}\pi)$ occurring in \eq{convective-} at least
for small pressures $|\pi|$ while allowing still for keeping positivity of
$K(\pi)$ for bigger negative pressures. Actually, dependence of $K$ on $\pi$
allows for modeling dependence of the sound speed on the pressure, which is
a well observed effect. More specifically, \eq{convective-} suggests the
expansion for small pressure as
$$
v_\text{\sc p}(\pi)=\sqrt{\frac{K{+}\pi}\varrho}
=\sqrt{\frac K\varrho}\sqrt{1+\frac\pi K}\sim
v_\text{\sc p,0}\Big(1+\frac\pi{2K}\Big)\ \ \ \text{ with }\ \ \
v_\text{\sc p,0}=\sqrt{\frac K\varrho}.
$$
Yet, actual fluids may exhibit even more pronounced dependence on pressure.
E.g.\ for water with $K=2.15\,$GPa and $v_\text{\sc p,0}=1500\,$m/s this
formula gives $v_\text{\sc p}(\pi)\sim
1500+0.348\times10^{-6}\pi$ while the measurements 
as in \cite[Tab.\,IV]{CheMil77SSSH} shows 
%
%
%
$v= 1521.46\,$m/s for $\pi=0$ and $v= 1538.12\,$m/s for $\pi=$100\,bar
=10\,MPa at temperature 20C, which gives $v(\pi)=v_0+ 2.146\times10^{-6}\pi$.
Even more pronounced dependence is in \cite[Tab.5]{NSBL16PPS},
reporting the compressibility at 0$^\circ$C in 10$^6$m$^3$/J=1/MPa is
$5.0792\times10^{-4}-1.3389\times10^{-6}\pi$ for $\pi$ in MPa.
 In general, sound speed provides high accuracy experimental
data exploited to for state equation from which the pressure dependence
of the bulk modulus (or compressibility) is determined, in the
particular case of  water cf.\ \cite{WagPru02IAPW} for the
Intl.\ Assoc.\ for the Properties of Water and Steam
1995 standard. 
\end{remark}






\begin{remark}[{\sl Uniqueness in the model \eq{Oskolkov-convective}}]\label{rem-uni}
  \upshape
  It is known that the incompressible Navier-Stokes equation enjoys
  uniqueness and full regularity if the pressure is in the class
  $L^1(I;L^{3/2}(\Omega))$, cf.\ \cite{Bers99SCRS} or
  also e.g.\ \cite{BerGal02RCIP}. Here, in the model \eq{Oskolkov-convective},
  the overall pressure $\frac1{2K}\pi^2+\pi$ is controlled even in
  $L^1(I;L^3(\Omega))$, and the additional force $\frac\varrho2({\rm div}\,\vv)\,\vv$ has a bilinear form similar to the convective term
  $\varrho(\vv{\cdot}\nabla)\vv$. There seems a chance to augment the
  uniqueness proof for the whole system \eq{Oskolkov-convective}
  involving also the other bilinear term $\frac1K\vv{\cdot}\nabla\pi$ in
  \eq{Oskolkov+}.
\end{remark}

\begin{remark}[{\sl
  Normal dispersion by convective conservative gradient terms}]
  \upshape
  To model dispersion in waves that can effectively propagate
  on long distances, one should use rather conservative gradient terms
  than dissipative ones in order not to make a substantial attenuation of
  the waves. For the normal dispersion, this can be realized through
  \eq{Sobolev}, cf.\ \eq{v=v(lambda)}. But the convective variant of
  \eq{Sobolev} seems problematic for analysis. Thus designing 
  a convective model which would realize normal dispersion through
  a conservative gradient term seems to be a challenging open problem.
\end{remark}

\begin{remark}[{\sl An alternative derivation}]
  \upshape
  Often, like \eq{NS-a} was written as \eq{variant-1}, 
  the system \eq{compress-NS} can equivalently be writen in a
  conservative
  instead of its convective form:
  \begin{subequations}\begin{align}\label{compress-NS1+}
    &\DT\rho+{\rm div}(\rho\vv)=0,
    \\&\label{compress-NS2+}
    \DT{\overline{\rho\vv}}+{\rm div}(\rho\vv{\otimes}\vv-\sigma)
    =\rho g\,.
  \end{align}\end{subequations}
  Then, the small perturbation ansatz \eq{ansatz} yields
\begin{align}\label{compress-NS2++}
  \varrho\Big(\!1{+}\frac\pi K\Big)\DT\vv+\varrho\frac{\DT\pi} K\vv
  +\varrho\frac{\nabla\pi}K{\cdot}(\vv{\otimes}\vv)
  +\varrho\Big(\!1{+}\frac\pi K\Big){\rm div}(\vv{\otimes}\vv)
  -{\rm div}\,\sigma=\varrho\Big(\!1{+}\frac\pi K\Big)g\,.
\end{align}
Substituting \eq{convective-}, the second term in \eq{compress-NS2++}
reads as $\varrho\frac{\DT\pi} K\vv=
-\varrho(1{+}\frac\pi K)({\rm div}\,\vv)\vv-\varrho\frac{\nabla\pi}K\cdot(\vv{\otimes}\vv)$
and, in particular, the tri-linear terms $\pm\varrho\frac{\nabla\pi}K\cdot(\vv{\otimes}\vv)$ cancel so that one obtains
\begin{align}\label{compress-NS2++-}
  \varrho\Big(\!1{+}\frac\pi K\Big)\DT\vv
  -\varrho\Big(1{+}\frac\pi K\Big)({\rm div}\,\vv)\vv
  +\varrho\Big(\!1{+}\frac\pi K\Big){\rm div}(\vv{\otimes}\vv)
  -{\rm div}\,\sigma=\varrho\Big(\!1{+}\frac\pi K\Big)g.
\end{align}
The assumption $|\pi|\ll K$ then leads to
\begin{align}\label{compress-NS2+++}
  \varrho\DT\vv-\varrho({\rm div}\,\vv)\vv
  +{\rm div}(\varrho\vv{\otimes}\vv-\sigma)=\varrho g\,.
\end{align}
By the calculus
  $(\vv{\cdot}\nabla)\vv={\rm div}(\vv\otimes\vv)-({\rm div}\,\vv)\vv$
  as in Remark~\ref{rem-conserv}, we obtain
 \begin{align}\label{compress-NS2++++}
   \varrho\DT\vv+\varrho(\vv{\cdot}\nabla)\vv-{\rm div}\,\sigma=\varrho g\,.
   \end{align}
  Using $\sigma=\bbD e(\vv)-\pi\bbI$ as before, this model however does not
  conserve the energy and, comparing it with
  \eq{NS-modif++}, we can see that the force
  $-\frac1K\pi\nabla\pi-\frac\varrho2({\rm div}\,\vv)\vv$ is to
  be added into the right-hand side of \eq{compress-NS2+++} to obtain a
  conservative model. To summarize, this alternative derivation
  leads to violation of energetics and needs to compensate it like
  before. 
It holds also for the diffusive variant \eq{Oskolkov-}, which
  would lead to enhancing the right-hand side of \eq{compress-NS2+} by
  $\frac KH(\Delta\rho)\vv$ and then
  \eq{compress-NS2++} and \eq{compress-NS2++-} would enhance by
  $\frac\varrho H(\Delta\pi)\vv$, which would again lead to
  \eq{compress-NS2+++} when using  \eq{convective-} enhanced  by the term
  $\frac1H\Delta\pi$ like \eq{Oskolkov+}.
\end{remark}

 \baselineskip=10pt

{\small
\section*{Acknowledgments}
The discussions with Donatella Donatelli, Eduard Feireisl, Alexander Linke,
Jos\'e-Francisco Rodrigues, Michael R{\accent23u}\v{z}i\v{c}ka, and
Giuseppe Tomassetti are warmly acknowledged.
}

\end{sloppypar}
\end{document}